\documentclass[a4paper,11pt,flqn]{article}
\usepackage{epsf}
\usepackage{enumitem}
\usepackage{graphicx}
\usepackage{mathptmx}
\usepackage{makeidx}
\usepackage{amsmath}
\usepackage{amssymb}
\usepackage{amsfonts}
\usepackage[bottom]{footmisc}
\usepackage{multicol}
\usepackage{type1cm}
\usepackage{mathrsfs}
\usepackage{arcs}
\usepackage{url}
\usepackage{color}
\usepackage[pagewise]{lineno}
\usepackage[francais]{minitoc}
\usepackage{setspace}
\usepackage{float}

\usepackage[T1]{fontenc}
\usepackage[applemac]{inputenc}
\newtheorem{theoreme}{Theorem}[section]
\newtheorem{lemme}[theoreme]{Lemma}
\newtheorem{proposition}[theoreme]{Proposition}

\newtheorem{definition}[theoreme]{Definition\rm}
\newtheorem{hypothese}[theoreme]{Hypothesis\rm}
\newtheorem{remarque}{\bf Remark}

\title{Entry-exit in the halo of a slow semi-stable curve}
\author{C. Lobry}
\date{\today}

\newcommand{\bitbul}{\begin{itemize}[label = \textbullet]}
\newcommand{\bittiret}{\begin{itemize}[label = -]}
\newcommand{\bito}{\begin{itemize}[label =$\circ$]}
\newcommand{\bit}{\begin{itemize}}
\newcommand{\fit}{\end{itemize}}
\newcommand{\ben}{\begin{enumerate}}
\newcommand{\fen}{\end{enumerate}}

\newcommand{\fin}{\end{document}}
\newcommand{\beq}{\begin{equation}}
\newcommand{\feq}{\end{equation}}
\newcommand{\dcom}{\begin{quote}\begin{small}}
\newcommand{\fcom}{\end{small}\end{quote}}

\newcommand{\bc}{\begin{center}}
\newcommand{\fc}{\end{center}}
\newcommand{\emat}{\mathrm{e}}
\newcommand{\ch}{\mathrm{cosh}}
\newcommand{\sh}{\mathrm{sinh}}
\newcommand{\eps}{\varepsilon}
\newcommand{\Rmat}{\mathbb{R}}
\newcommand{\Nmat}{\mathbb{N}}

\newcommand{\arc}{\overset{\displaystyle  \frown}}
\def\1{{\rm 1\mskip-4.4mu l}}
\newcommand{\eset}[1]{{
  \mathchoice
    {\left\{\!\!\left\{ #1 \right\}\!\!\right\}} 
    {\left\{\!\left\{ #1 \right\}\!\right\}} 
    {} 
    {} 
  }
}

\begin{document}
\maketitle

\texttt{ This text is the translation (with the help of DeepL) of the article {\em Entrée-sortie dans le halo d’une courbe lente semi-stable} of the same author  posted at \url{
https://doi.org/10.48550/arXiv.2203.04712
}}

\section*{Introduction} 

I am interested in the phase portrait of the system:
 \beq \label{Smintro}
 S_m\quad  \quad \left\{
\begin{array}{lcl}
\displaystyle \frac{dx}{dt}& =& 1\\[6pt]
\displaystyle  \frac{dy}{dt} &=& \displaystyle \frac{1}{\eps}\sqrt{m^2+y^2} \Big( f(x) -y\Big)
 \end{array} 
 \right.
\feq
when the parameters $m$ and $\eps$ are small. This system appears in a natural way in the study of the important  phenomenon called {\em inflation} in population dynamics (see \cite{BLSS21,HOLTPNAS20,KAT21}), which motivates my study.

Some forty years ago, particular solutions of systems of the form :
 \beq \label{eqLR2}
 \mathrm{LR}\quad \quad  \quad \left\{
\begin{array}{lcl}
\displaystyle \frac{dx}{dt}& =& F(m,x,y)\\[6pt]
\displaystyle  \frac{dy}{dt} &=& \displaystyle \frac{1}{\eps}G(m,x,y)
 \end{array} 
 \right.
\feq
\begin{minipage}{0.52 \textwidth}
were defined, in the particular case where the set  $G(0,x,y) = 0$ consists of two curves $\gamma_i$, called ''slow curves'', which intersect transversally and where $G$ changes sign at the crossing of each $\gamma_i$. In this case each curve $\gamma_i$  decomposes into $\gamma_i^{\;a}$ (attractive part) and $\gamma_i^{\;r}$ (repulsive part). 
 When $\eps$ is small one shows that for some $m$ of the order of $\emat^{-\frac{1}{\eps}}$ there exist, surprisingly, some 
 \end{minipage} $\quad$ \begin{minipage}{0.45 \textwidth}
\includegraphics[width=1\textwidth]{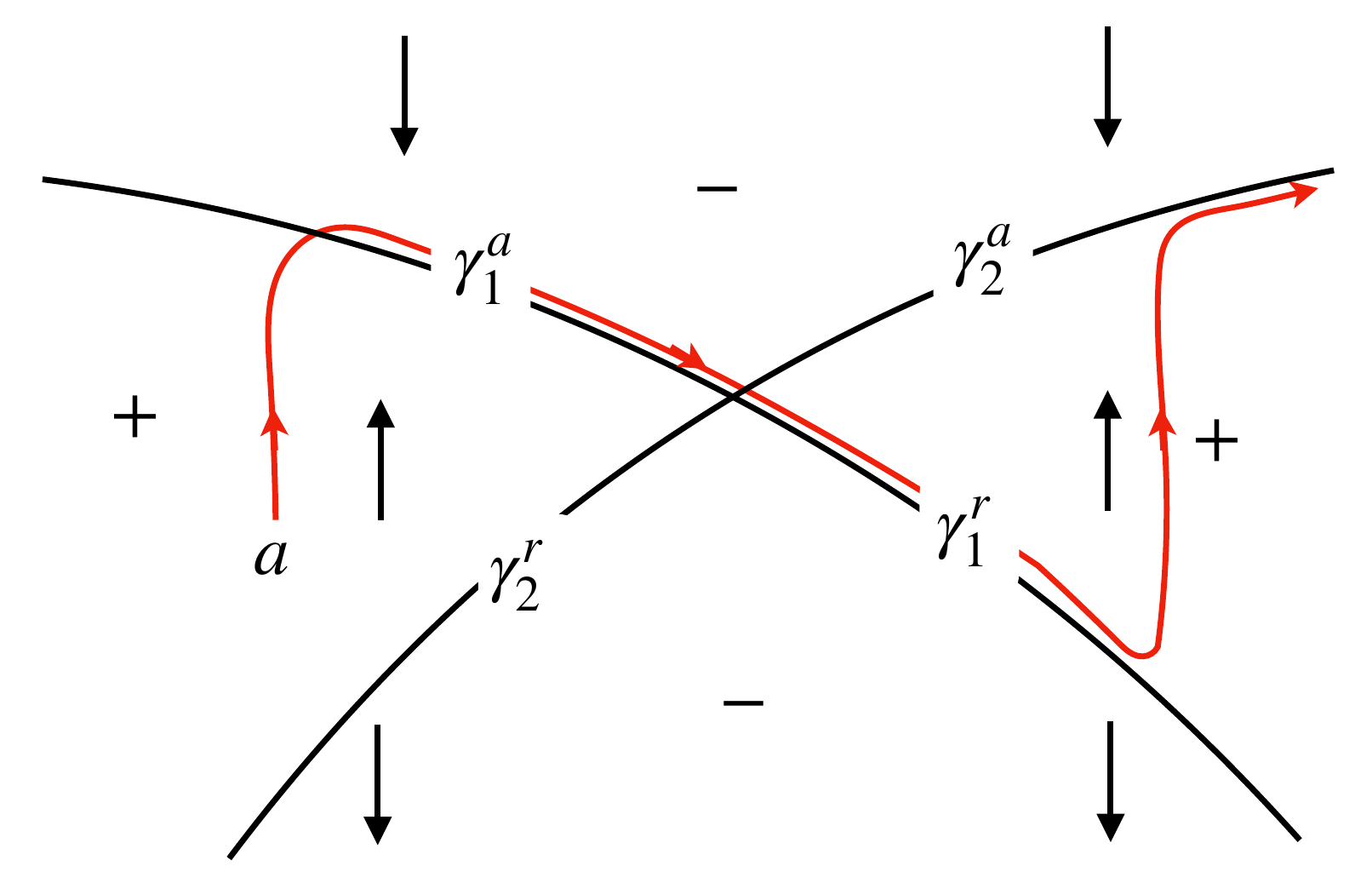}
\end{minipage} \\[2pt]
  solutions which, like the solution in the figure above, after having followed the attractive part $\gamma_1^{\;a}$ of the slow curve $\gamma_1$ continue to follow, for a significant time, the repulsive part $\gamma_1^{\;r}$ (see the figure). These solutions have been called {\em canards} by their inventors: (see \cite{BEN81, BCDD81}). One can find online on Scholarpedia the article \cite{WEC07} which gives the history of this discovery and its developments. The interested reader will also find in \cite{FRU09} a review of the applications of the {\em canards} theory to the study of dynamic bifurcations.\\

The system I am interested in is similar to the one I have just discussed but not exactly the same. For $m=0$ we have :
 \beq 
 S_0\quad  \quad \left\{
\begin{array}{lcl}
\displaystyle \frac{dx}{dt}& =& 1\\[6pt]
\displaystyle  \frac{dy}{dt} &=& \displaystyle \frac{1}{\eps}|y| \Big( f(x) -y\Big)
 \end{array} 
 \right.
\feq
Here, it is $ \displaystyle \frac{1}{\eps}|y| \Big( f(x) -y\Big)$ which corresponds to the function $G$, the slow curves $\gamma_1$ and $\gamma_2$ are respectively the graph of $f$ and the axis $y = 0$ but the distribution of  \\[2pt]
\begin{minipage}{0.50 \textwidth}
signs are different.
 Thus the slow curve $\gamma_1$ does not go from attractive to repulsive as before but remains attractive whereas $\gamma_2$ (here the axis $y = 0$ ) is semi attractive, semi repulsive. I am  are going to show the existence, always for exponentially small values of $m$ with respect to $\eps$, of solutions which, like the one in the figure, run along portions of $\gamma_2$, that is to say portions {\em non stable} of the slow curve,
\end{minipage}$\quad $
\begin{minipage}{0.48 \textwidth}
\includegraphics[width=1\textwidth]{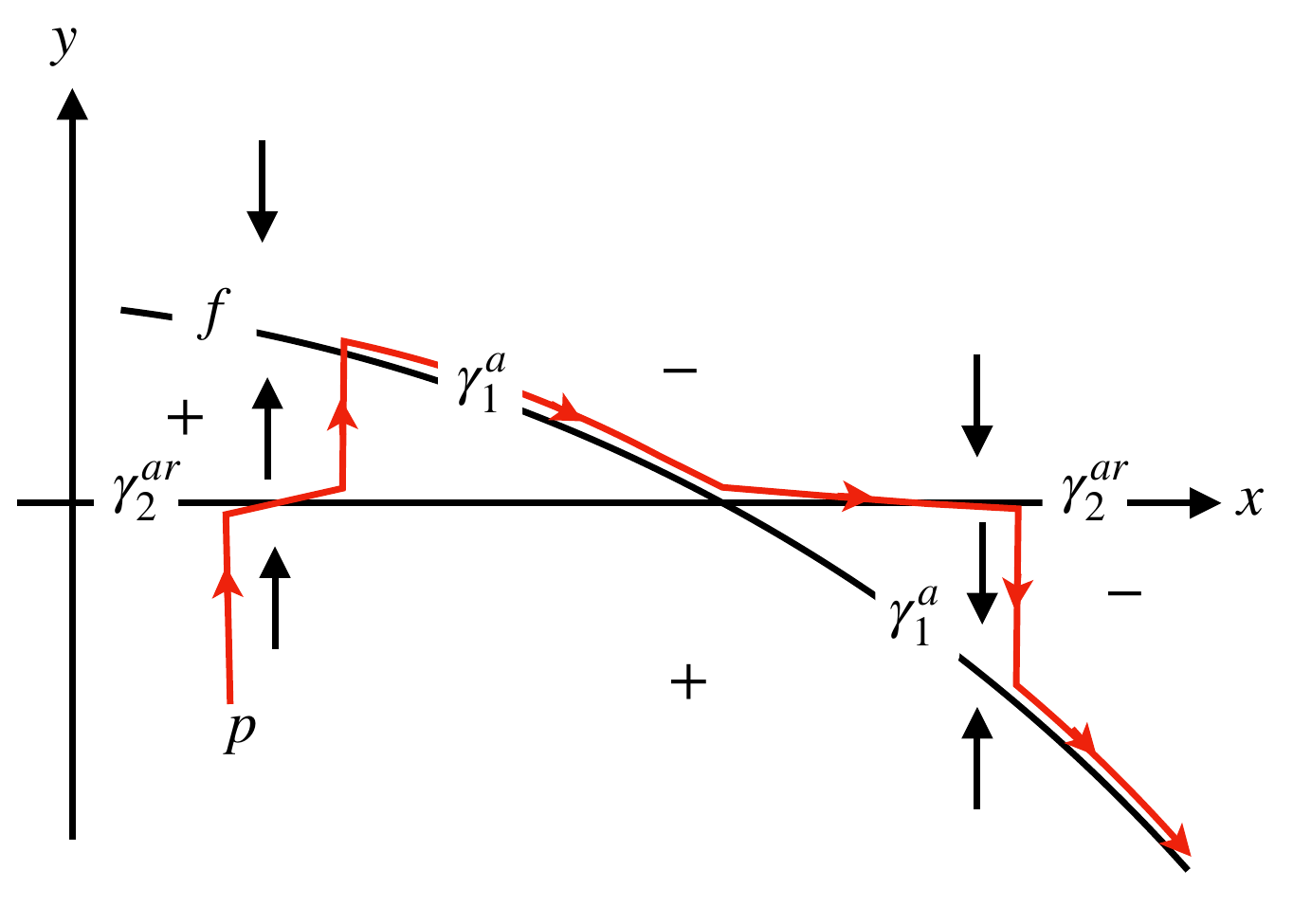}
\end{minipage}\\[6pt]
  {\em for significant durations}. We can no longer speak of  ''canard'' in the sense of the usual definition but the phenomenon is of the same nature as I will show.

 To prove the existence of these solutions I will use the methods of the inventors of the ''canards'', i.e. Non-Standard Analysis (NSA) in Nelson's I.S.T. formalism and the techniques of changes of variables developed in particular in \cite{BEN81, BCDD81}. 
 The reader familiar with these papers will notice that the present article is a simple application of the methods introduced there. Nevertheless, to make this text understandable to a wider audienceI have tried to keep the concepts and terminology related to NSA to a minimum and have described them in a short appendix.  
 
  The objective is to explain the phase portrait of \eqref{Smintro} and, because of the possible applications, already mentioned, in population dynamics, I deal with the case where $f$ is not necessarily continuous.     After having specified the hypotheses and the notations in the second section I state and illustrate on examples the main result of this article: the theorem of approximation of solutions.   In the next section I prove the theorem and then I finish with an application : the proof of a conjecture stated by G. Katriel in \cite{KAT21}. The demonstrations are geometrical and are essentially done on figures which, if the article is read on a double page, are as much as possible opposite to the text. 
\newpage  \tableofcontents
\newpage
\section{The $S_m$ system}
\subsection{Vocabulary and notations} 
It is necessary to introduce right away a minimal NSA vocabulary and the associated notations. The reader who does not know the NSA methods can be satisfied  with the intuitive meaning of the introduced words but he will find, if he wishes, more complete information in the appendix.

In non-standard analysis, all is about the fact  that the sentence : 
 \bitbul
 \item $\eps$ is a strictly positive real number infinitely small (i.e. {\em infinitesimal} ) 
 \fit
 has a precise formal meaning (see appendix \ref{ANS}). 
 
 From this we immediately define :
 \bitbul
  \item $x \in \Rmat$ is {\em infinitely large} $\stackrel{df}{\longleftrightarrow}  \exists \, \eps\; \mathrm{infinitesimal\; such\;that}\; x > \frac{1}{\eps}$ 
 \item $x \in \Rmat$ is {\em limited} if $\stackrel{df}{\longleftrightarrow}  |x|$ is not infinitely large..
\fit
 and we introduce the notations :
 \bitbul
 \item $x \sim y \stackrel{df}{\longleftrightarrow} |x-y|$ is infinitesimal.
 
 \item $x \lnsim y \stackrel{df}{\longleftrightarrow} x<y$ and $ |y-x] \nsim 0$
 
\item The {\em halo} of a subset $E$ of $\Rmat^2$ consists of the points which are infinitely close to the points of $E$, i.e. $d(x,E) \sim 0$.

\item The infinitely large and infinitely small reals are {\em nonstandard} but there are many others such as, for example, $1+\eps$ when $\eps$ is infinitely small\footnote{
Numbers in decimal form with a finite number of decimal places provide a good picture: we can decide that numbers with 6 decimal places are standard and that numbers with 12 decimal places are nonstandard. An infinitesimal is a number of the form $0.000;000;***;***$.}. Every limited real number is infinitely close to a unique standard number.

\item A function $f$ is {\em limited} if there exists $M$ limited such that  $\forall x \;|f(x)|< M$ and  $C^1$-{\em limited } if it, and its derivative, are both limited.
 \fit
 
\subsection{The equations}\label{equationsdebase}
In this article we focus on the system  :
 \beq \label{Sm}
 S_m\quad  \quad \left\{
\begin{array}{lcl}
\displaystyle \frac{dx}{dt}& =& 1\\[6pt]
\displaystyle  \frac{dy}{dt} &=& \displaystyle \frac{1}{\eps}\sqrt{m^2+y^2} \Big( f(x) -y\Big)
 \end{array} 
 \right.
\feq
\paragraph{Hypotheses and notations on $f$:}\label{hyp}
We assume that $f$ is piecewise $C^1$-limited, i.e. :
\ben
\item There exists a discrete sequence  $D = \{x_n \;;\; n\in \mathbb{Z^*}\}$ 
\item On each  interval $[x_n,x_{n+1}[$ the function  $f$ is the restriction of a function   $\tilde{f}$, $C^1$-{\em limited,} defined on the whole  $\Rmat$ such that  $\tilde{f}(x) = 0 \Rightarrow f'(x) ≠0$ (by the way the zero's of  $\tilde{f}$ are isolated).
\item  Let $\mathcal{C}$ be the graph of $f$ and $\theta_n,:\, n \in J$ ($\theta_n < \theta_{n+1}$) the discrete sequence  of the values where  $f$ changes  sign. For a given $x$ we note $\theta(x)$ the first $\theta_n$ larger than $x$.
\fen
 
  \begin{definition}\label{c-graph} The graph of $f$ with the addition of the vertical segments that join the left and right limits of $f$ at a point of discontinuity is called the complete graph of $f$, or  the C-graph.
\end{definition}

Under these assumptions\footnote{
These hypotheses are not absolutely necessary, but they allow us to simplify the writings, in particular by implying that the solutions cannot tend to infinity for finite values of $t$.}, 
for any initial condition, the system $S$ admits a unique solution defined for all $t$. Precisely if $x_0 \in [\theta_n,theta_{n+1}[$ we integrate $S_m$, with $\tilde{f}_n$ instead of $f$, from the initial condition $(x_0,y_0)$ ; as $\tilde{f}_n$ is bounded, the solution is defined until it meets the vertical $x = x_{n+1}$ at a point $(x_{n+1}, y_{n+1})$ which serves as a new initial condition and so on.
\bitbul
\item Let $(x(t, x_0,y_0,m),y(t,x_0,y_0,m))$ be the solution of $S_m$ of initial condition $(x_0,y_0)$ at time $0$. We have $x(t, x_0,y_0,m) = x_0+t$.
\fit
Note that this existence of solutions theorem is true for any $\eps >0$ and any $m$, whether they are standard or not. From now on:
 \bitbul
\item It is assumed that $\eps >0$ { infinitely small} is given once and for all and 
$$m= \emat^{\frac{\rho}{\eps}\quad \quad \rho < 0}$$
 is a parameter. 
\fit
We are interested in the evolution of the phase portrait of $S_m$ as a function of $\rho$.

Note that it is not the concern for maximum generality that leads us to consider piecewise continuous-functions rather than simply continuous ones, but contexts where modeling with discontinuous functions is more natural (see \cite{BLSS21}). Moreover this greater generality does not introduce any additional difficulty in the type of proofs I use.

\subsection{The constrained system and the approximation theorem}\label{contraint}
\begin{figure}[t]
  \begin{center}
 \includegraphics[width=0.9\textwidth]{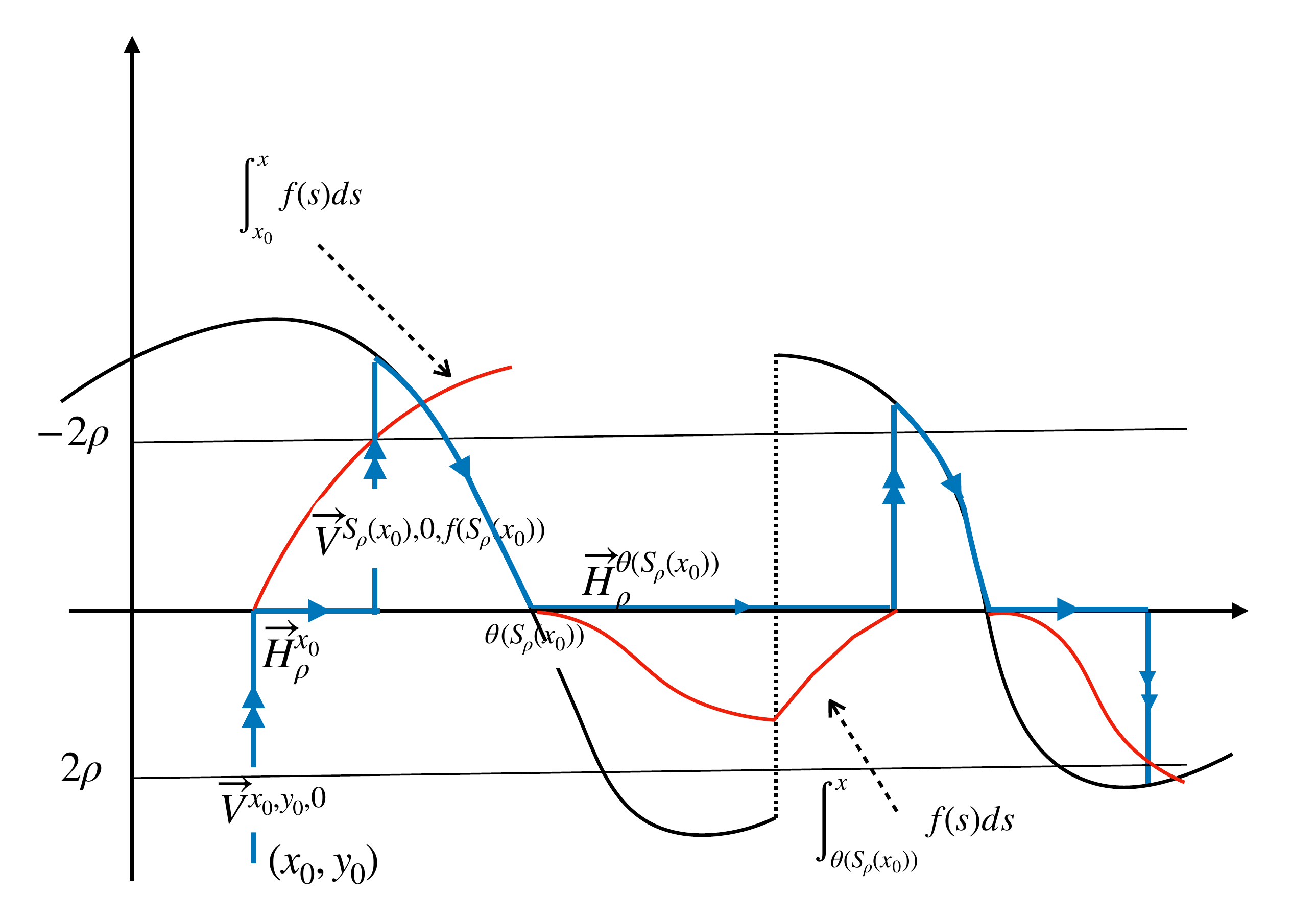}
 \caption{C-trajectory of \eqref{SCm} from $(x_0,y_0)$. ; explanations subsection \ref{contraint} } \label{figSCm}
 \end{center}
 \end{figure}

To the system $S_m$ we associate the ''constrained system'' :
 \beq \label{SCm}
 S_m^0\quad \quad \quad \left\{
\begin{array}{lcl}
\displaystyle \frac{dx}{dt}& =& 1\\[6pt]
\displaystyle 0 &=& \displaystyle \sqrt{m^2+y^2} \Big( f(x) -y\Big)
 \end{array} 
 \right.
\feq
It is not a differential system but, however, one can associate ''pseudo-trajectories'' to it in the following way\footnote{
The definition of the notion of solution for "constrained systems", of the form above, but in higher dimension, has been the subject of important studies (for example \cite{TAK76}) but it is not necessary to refer to them here because they do not deal with the link between the solutions of $S_m$ and $S^0_m$. On the other hand one can consult \cite{WEC07} and its references on the existence of {\em ducks} in dimension greater than 2.}.
\begin{definition} We call :
\bitbul
\item \textbf{Vertical segment} : a segment of the form $$ \overrightarrow{V}^{x,[a,b]}=[(x,a), (x,b)]$$
such that $0 \not \in ]a,b[$ and $b = 0$ or $b = f(a)$ oriented upward below the graph of $f$ and downward above. 
	
\item \textbf{Slow curve segment} : a part of the graph of $f$ of the form:
$$ \mathcal{C}^x = \{(s,f(s))\,:\, x \leq s \leq \theta (x)\}$$ 
where $\theta(x)$ is the first change of sign of $f$ that follows $x$ (cf. assumptions \ref{hyp}.3.), oriented to the right.

\item \textbf{Horizontal segment} : a rightward-oriented segment of the form $\overrightarrow{H}_{\rho}^{x} = [(x, 0),(S_{\rho}(x),0)]$ where $S_{\rho}(x)$ is defined as the smallest (possibly $+\infty$) $x^*>x$ such that :
\beq
\displaystyle \int_x^{x^*} f(s)ds \in \{+2\rho,\,0, -2\rho \}
\feq
when $\rho < 0$. The end $(S_{\rho}(x),0)$ of the horizontal segment is called the "exit point" (from the halo of $\mathcal{C}$).
\fit
\end{definition}
\begin{proposition} \label{rderho} 
For fixed $x$, $S\mathcal{C}(x)$ is a decreasing function of $\rho$ which tends to $x$ when $\rho$ tends to $0$.
\end{proposition}
\textbf{Proof.} Obvious. $\Box$
\begin{definition}
A \textbf{C-trajectory} (for trajectory of the " constrained system ") of (\ref{SCm}) is a sequence $\Gamma_i$ of segments put together (i.e. the origin of $\Gamma_{i+1}$ is the end of $\Gamma_i$) such that:
\bito

\item to a vertical segment $ \overrightarrow{V}^{x,[a,f(x)]}$ succeeds the slow curve segment $ \mathcal{C}^x$ 
\item to a vertical segment $ \overrightarrow{V}^{x,[a,0]}$ succeeds the horizontal segment $\overrightarrow{H}_{\rho}^{x}$
\item to a horizontal segment $\overrightarrow{H}_{\rho}^{x}$ succeeds the vertical segment $ \overrightarrow{V}^{S_{\rho}(x),[0,f(S_{\rho}(x))]}$
\item a slow curve segment $ \mathcal{C}^x$ is followed by the horizontal segment $\overrightarrow{H}_{\rho}^{\theta(x)}$
\fit
These rules allow to associate to any ''initial condition'' $(x_0,y_0)$ a unique C-trajectory as can be seen on figure \ref{figSCm}.
\end{definition} 
\begin{definition}\label{longe} Let $\Gamma$ be an arc $s \mapsto (\alpha (s),\beta(s))$, ($s \in [a,b]$), of $\Rmat^2$ and $t \mapsto (x(t),y(t))$, ($ t\in [t_1,t_2]$), another arc of $\Rmat^2$. We say that 
$ (x(t),y(t)) \;\mathrm{hugs}\; \Gamma$
if there exists a parameterization $t \mapsto s(t)$, ($t \in [t_1,t_2]$), of $\Gamma$ such that, 
for all $t$, $d((x(t),y(t)),(\alpha(s(t),\beta(s(t)))) \sim 0$, where $d$ is the natural distance of $\Rmat^2$.
\end{definition}
We can now state the main result of this paper.

\begin{theoreme} \label{theorem} \textbf{Approximation theorem.}
Let $(x_0,y_0)$ be a limited initial condition, such that $y_0 \not \sim 0$, and $t \mapsto (x(t),y(t))$ be the solution of the differential system \eqref{Sm} issued from it at time $t_0$. Then $(x(t),y(t))$ is infinitely close to the C-trajectory of \eqref{SCm} coming from the same point. More precisely, if $\,\Gamma_1, \Gamma_2,\cdots, \Gamma_n,\cdots$ are the successive segments of the C-trajectory, there exists a sequence of instants, $t_1,t_2,\cdots,t_n,\cdots$ such that on $[t_{n-1},t_n]$ the trajectory $(x(t),y(t))$ hugs (see def. \ref{longe}) the segment $\Gamma_n$.
\end{theoreme}
The denomination ''exit point'' for the extremity $(S_{\rho}(x),0)$ of the horizontal segment $\overrightarrow{H}_{\rho}^{x} $ comes from the fact that the trajectory which  ''entered'' at the point $(x,0)$ in the halo of $y = 0$ ''exit out'' at the point $(S_{\rho}(x),0)$. For this reason we call the quantity $R_{\rho}(x) = S_{\rho}(x) - x$ "delay" (understood " for the exit").

\textbf{The hypothesis $y_0 \not \sim 0$ is essential as we shall see in the paragraph \ref{ci0}}

\paragraph{Note} (For those who know I.S.T.) Using all the power of I.S.T. we can introduce the concept of {\em shadow} of a set (the {\em standardised} of the {\em  halo}). The previous speech is the rephrasing of the proposition: {\em The shadow of the trajectory issued from the point $a$ is the C-trajectory issued from $a$}. By limiting ourselves to a simplified version of I.S.T. (as I do here) one deprives oneself of tools facilitating the writing as would be the use of the shadow but decreases the price to pay for the practice of NSA. I discuss this point in the appendix.

\subsection{Illustrations of the approximation theorem}
\begin{figure}
  \begin{center}
 \includegraphics[width=0.9\textwidth]{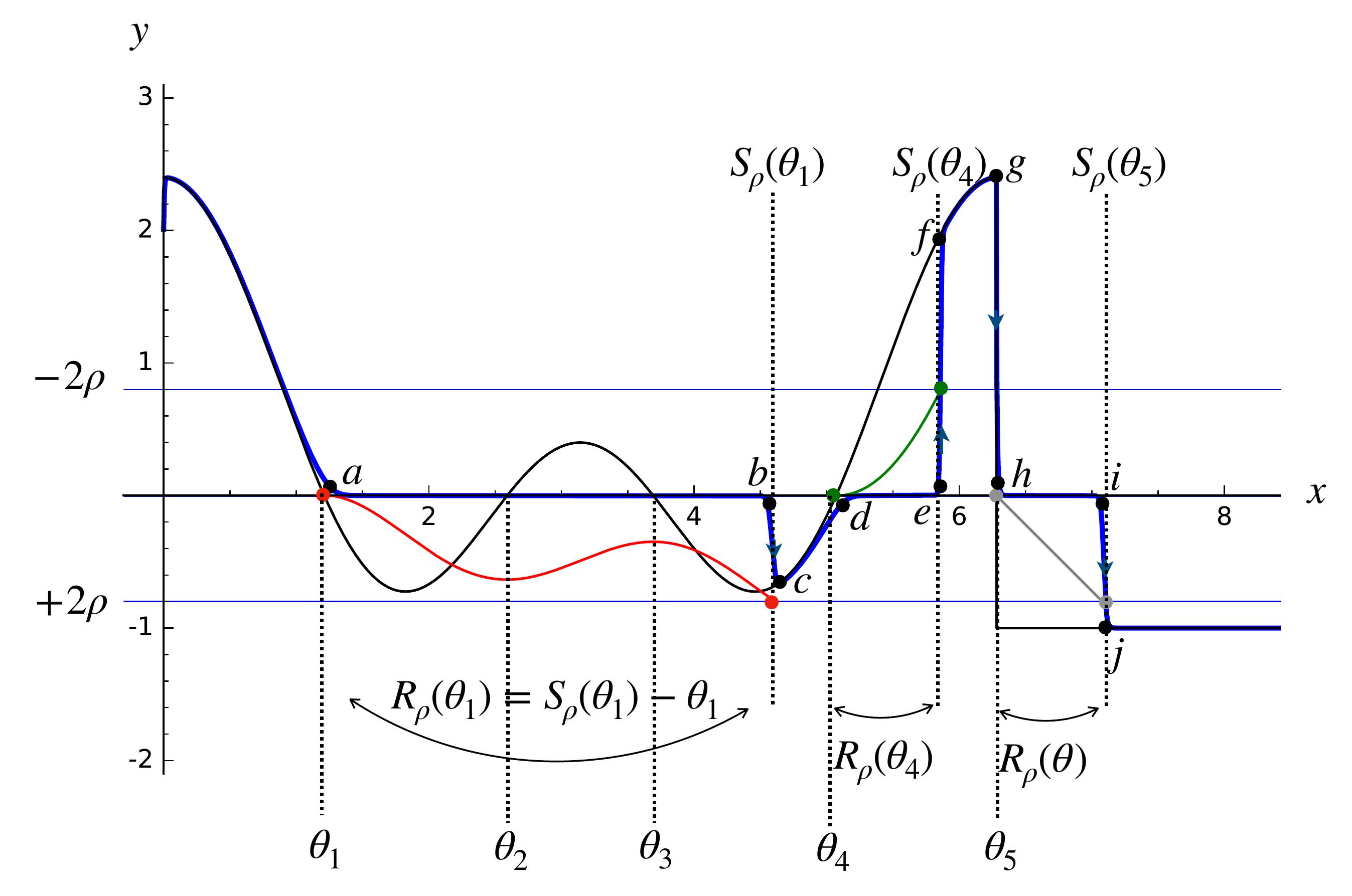}
 \caption{Explanations in text} \label{retard5}
 \end{center}
 \end{figure}

 \paragraph{Comments on figure \ref{retard5}. }
One sees a simulation of the system  $S_m$  with :
 \bitbul
 \item The function $f$ defined by :
 \beq
 \left\{
\begin{array}{lcl}
\displaystyle x \leq 2\pi& \Rightarrow & f(x) = \cos(x)+\cos(2x)+0.4 \\
\displaystyle x> 2\pi& \Rightarrow & f(x) = -1
\end{array}
\right.
\feq
Its graph is drawn in black ; $f$ is null for the values $\theta_1,\theta_2,\theta_3$ and $\theta_4$ which it is useless to specify and it is discontinuous for $x = 2\pi$.
\item $\eps = 0.01\quad \rho = -0.4\quad \displaystyle m = \emat^{-\frac{0.4}{\eps}}=4.24\, 10^{-18}$
 \fit
 In thick blue I have simulated the trajectory $(x(t),y(t))$ coming from the point $(0,2)$ at the time $t_0 = 0$. We can see that, apart from the "slightly rounded angles", it presents a sequence of "segments" very close to the C-trajectory. We detail this sequence :
 \ben

 \item There exists $t_1 \sim 0$ such that on $[t_0,t_1]$ the solution hugs the vertical segment $\overrightarrow{V}^{0,[2, 2.4]}$.
 Let $\theta_1$ be the first value after $x = 0$ for which $f$ changes sign ; there exists $t_2 \sim \theta_1$ such that on $[t_1,t_2]$ the solution hugs the slow curve segment $\mathcal{C}^0 = \{ (x,f(x)\;;\; x \in[0,\theta_1]\}$ up to point $a \sim(\theta_1,0)$. The point $a$ is an entry point into the halo of the line $y = 0$.

 \item Now we determine the exit point following $(\theta_1,0)$. In red we have drawn the graph of the function:
 $$x \mapsto \int_{\theta_1}^xf(s)ds$$
 which meets the line $y = 2\rho$ at the point of abscissa $S_{\rho} (\theta_1) = \theta_1+R_{\rho} (\theta_1)$ ; There exists $t_3 \sim S_{\rho} (\theta_1)$ such that on $[t_2,t_3]$ the trajectory hugs the horizontal segment $[(\theta_1,0),(S_{\rho} (\theta_1),0)]$ up to the point $b \sim (S_{\rho} (\theta_1),0)$. Note that there is no reason for $S_{\rho}(\theta_1)$ to precede $\theta_2$ or $\theta_3$.

 \item There exists $t_4 \sim S_{\rho} (\theta_1)$ such that on $[t_3,t_4]$ the solution hugs the vertical segment $\overrightarrow{V}^{S_{\rho} (\theta_1),[0,f(S_{\rho} (\theta_1))]}$ up to the point $c \sim (S_{\rho} (\theta_1), f(S_{\rho} (\theta_1))$.
 Let $\theta_4$ be the first value after $x = S_{\rho} (\theta_1))$ for which $f$ changes sign; there exists $t_5 \sim \theta_4$ such that on $[t_4,t_5]$ the solution hugs  the slow curve segment $\mathcal{C}^{S_{\rho} (\theta_1)} = \{ (x,f(x)\;;\; x \in[S_{\rho} (\theta_1),\theta_4]$ up to point $d \sim(\theta_4,0)$. The point $d$ is an entry point in the halo of the line $y = 0$.

 \item Determination of the exit point following $S_{\rho} (\theta_4,0)$. In green we have drawn the graph of :
 $$x \mapsto \int_{\theta_4}^xf(s)ds$$
 whose intersection with the line $y = - 2\rho$ determines the exit point $(S_{\rho} (\theta_4),0)$. There exists $t_6\sim S(\theta_4)$ such that on $[t_5,t_6]$ the solution hugs the horizontal segment $H_{\rho} ^{\theta_4}$.
 We leave it to the reader the exercise  of interpreting  the following trajectory segments: $ef,\; fg,\; gh,\; hi,\; ij$
 \fen 

\noindent \textbf{Remark}. 
 The predictions of the approximation theorem are very precise whereas the quantities manipulated in the numerical scheme are remarkably small as for example $m = 4.24\, 10^{-18}$ ; this is because the representation of the reals in "floating point" notation allows to represent very small numbers ($10^{-250}$ with my software) If we make the change of variable $z = y-1$ which translates the figure around the line $z = 1$ the predictions of the approximation theorem would be much less efficient ; this point is detailed in \cite{LOB92}.

\begin{figure}
  \begin{center}
 \includegraphics[width=0.8\textwidth]{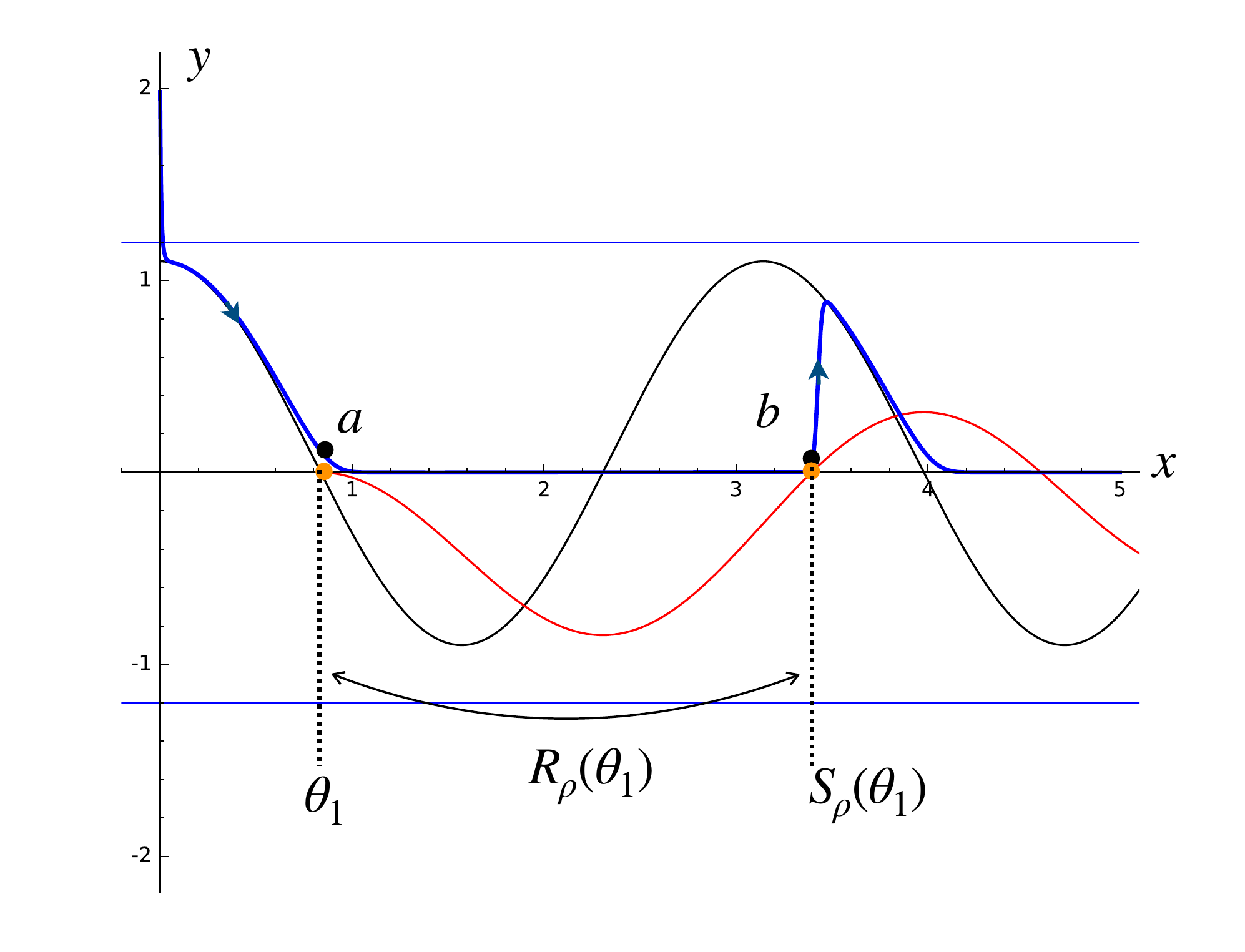}
 \caption{Explanations in text} \label{retard10}
 \end{center}
 \end{figure}

\newpage
  \paragraph{Comments on the figure \ref{retard10}. } 
 On the example of the figure \ref{retard5}, when the trajectory hugs a horizontal segment it always crosses the axis $y=0$. This is not the case when the exit point is determined by :
 $$ \int_x^{S_{\rho} (x)}f(s)ds = 0$$
 On this example we took $\eps = 0.01$, $f(x) = 0.5\cos(x) +0.1$ which cancels a first time (for $x \geq 0$) in $\theta_1 = \arccos (-0. 2) $ and we took $\rho = -0.6$ so that the graph (in red) of $x \mapsto \int_{\theta_1}^xf(s) ds$ intersects the axis before meeting $y = ±2\rho$.
 
\begin{figure}[h!]
  \begin{center}
 \includegraphics[width=0.9\textwidth]{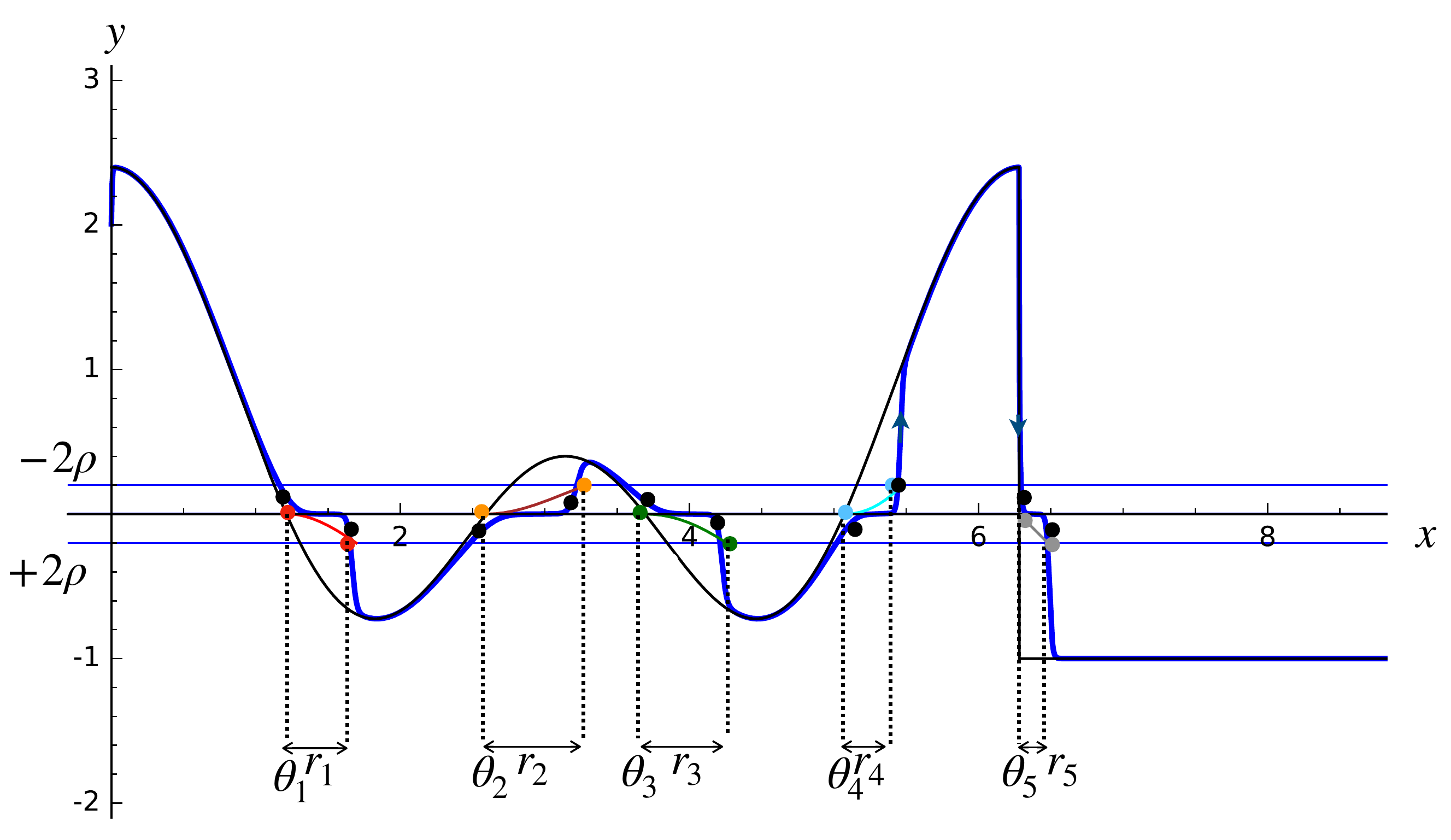}
 \caption{Explanations in text} \label{retard5bis}
 \end{center}
 \end{figure}

 \paragraph{Comments on figure \ref{retard5bis}. } 
  For this simulation we used the model of the figure \ref{retard5} except $\rho = - 0.1$ instead of $\rho = -0.4$. We can see that the solution is qualitatively quite different: it now presents 5 horizontal segments, instead of 3 in the case $\rho = -0.4$. The reader will easily convince himself that, in general, the sum of the lengths of the horizontal segments tends to $0$ when $\rho$ tends to $0$.

\begin{figure}
  \begin{center}
 \includegraphics[width=1\textwidth]{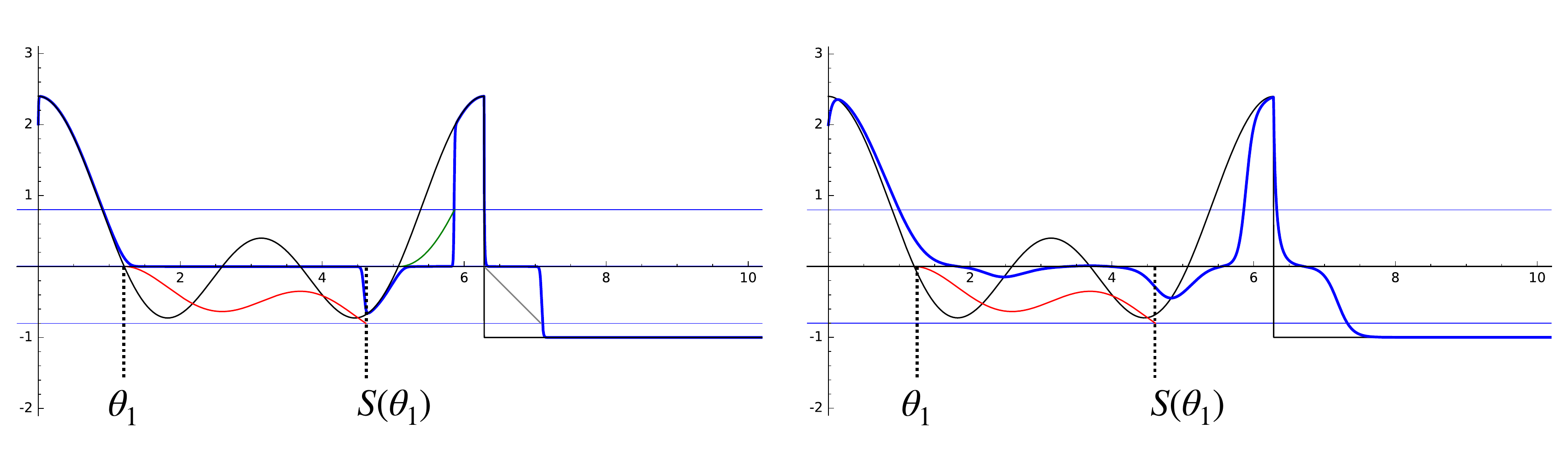}
 \caption{Left $\eps = 0.01$, right $\eps = 0.1$} \label{retard6bis}
 \end{center}
 \end{figure}

 \paragraph{Comment on the figure  \ref{retard6bis}.}
  
The approximation theorem is stated with $\eps$ infinitely small. For the simulations of the figure \ref{retard5} we took $\eps = 0.01$. What about a ''bigger epsilon'' ? On this figure we compare the simulation obtained for $\eps = 0.01$ with what the simulation of the same system gives with $\eps = 0.1$. The successions of horizontal and vertical segments are less clear but the predictions of the output points remain valid.

\section{Proof of the approximation theorem}
This section is devoted to the proof of the approximation theorem. It consists in the examination of phase portraits and the following of typical trajectories. We will leave it to the reader to convince himself that all possible cases have been considered.

\subsection{Outside the halo of the horizontal axis $y = 0$.}
In the figure on the right we have represented the graph of $f$ (in black) in a zone where it does not meet the axis $y = 0$. We consider two initial conditions, $(x_1,y_1)$ below the graph, not infinitely close to the graph of $f(x)$ nor to the axis $y = 0$ and $(x_2,y_2)$ above the graph.\\
\begin{minipage}{0.49 \textwidth}
 The green curves are the respective graphs of $f(x)+\alpha$ and $f(x)-\alpha$ where $\alpha$ is chosen strictly positive \textbf{non infinitely small}. The segments $[a,b]$ and $[a,c]$ are respectively carried by the lines of slope $-1/\alpha$ and $+ 1/\alpha$ coming from the point $(x_1,y_1)$ until they meet the graph of $f(x)-\alpha$. We thus define a domain $D$ bounded by the 
\end{minipage} $\quad$
\begin{minipage}{0.5 \textwidth}
\includegraphics[width=1\textwidth]{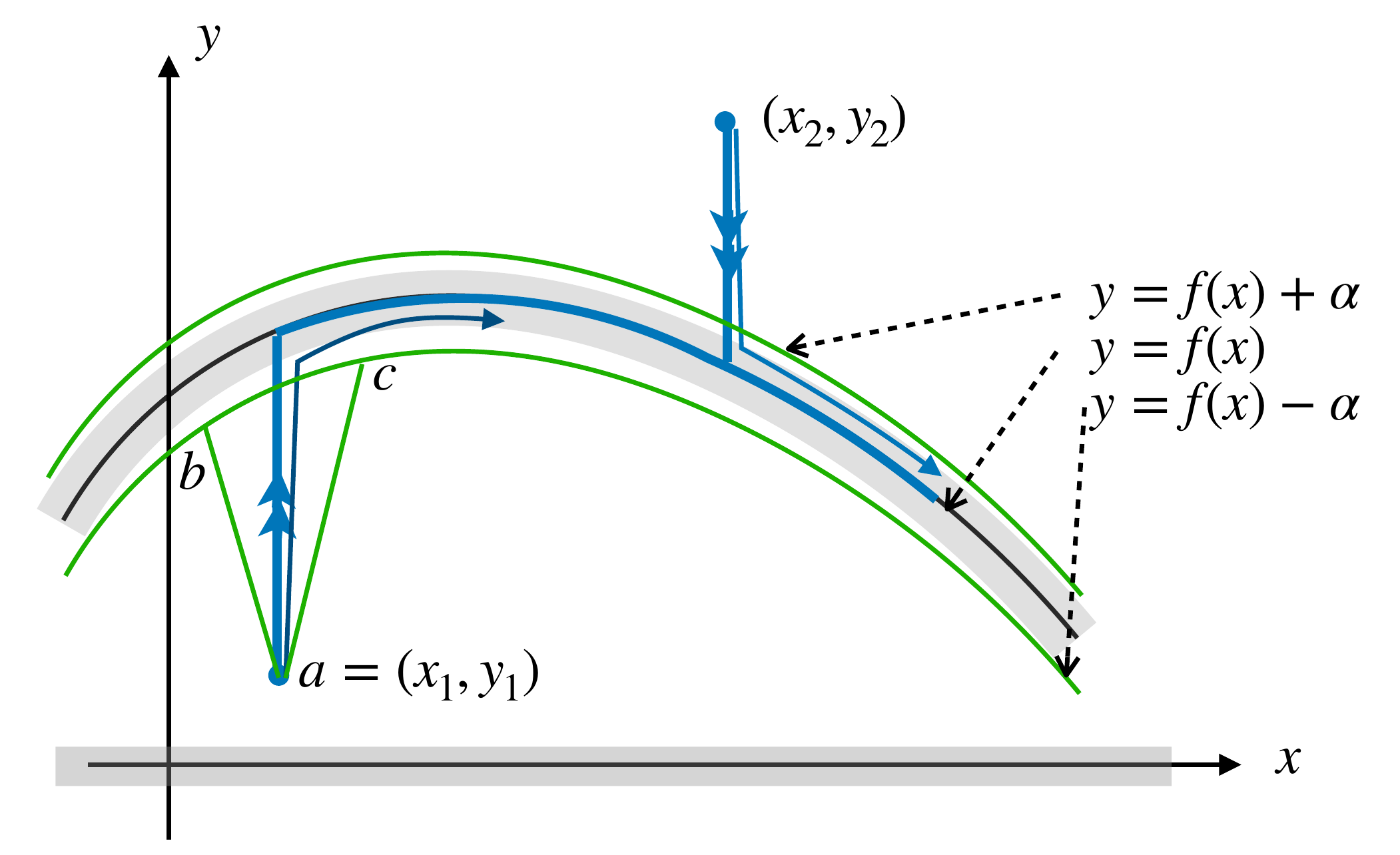}
\end{minipage} 
segments $[a,b]$, $[a,c]$ and the arc $(a,b)$ carried by the graph of $f(x)-\alpha$. Consider the system $S_m$ :
 \beq \label{Sm1}
 S_m\quad \quad \left\{
\begin{array}{lcl}
\displaystyle \frac{dx}{dt}& =& 1\\[6pt]
\displaystyle \frac{dy}{dt} &=& \displaystyle \frac{1}{\eps}\sqrt{m^2+y^2} \Big( f(x) -y \Big)
 \end{array} 
 \right.
\feq
Since the domain $D$ does not encounter a grey area, neither the factor $\sqrt{m^2+y^2} $ nor the factor $(f(x)-y)$ is infinitely small and therefore the second member of the second equation of \eqref{Sm1} is infinitely large positive and hence along the segments $[a, b]$ and $[a,c]$ the field points strictly inside the domain $D$ and, as it does not cancel, the trajectory coming from $a = (x_1,y_1)$ leaves  $D$, after a time necessarily infinitely small, at a point of the arc $(a,b)$. As this holds for all $\alpha \gnsim 0$ it is clear that the trajectory will end up being infinitely close to the graph of $f(x)$, on the other hand it is less clear that there exists a $t_1$ \textbf{infinitely small} for which $(x(t_1),y(t_1))$ is infinitely close to the graph of $f(x)$ because as one gets closer to the graph the speed of $y$ decreases.
This point is acquired by a typically non-standard argument detailed in the appendix \ref{halocourbelente}.
It is easy to convince oneself that once the solution $(x(t),y(t))$ has entered the halo of the graph of $f(x)$, the trajectory must remain there as long as $x(t)$ remains non-infinitely close to a value where $f(x)$ changes sign.

\subsection{In the halo of the horizontal axis $y = 0$ : the exponential magnifying glass}
\begin{definition}
We call {\em the exponential magnifier} the change of variable:
$$z = [y]^{\eps} \;\stackrel{df}{=} \;\mathrm{sgn}(y)|y|^{\eps}$$
\end{definition}
This change of variable was introduced in \cite{BEN81} and \cite{BCDD81} .\\[4pt]
\begin{figure}
  \begin{center}
 \includegraphics[width=1\textwidth]{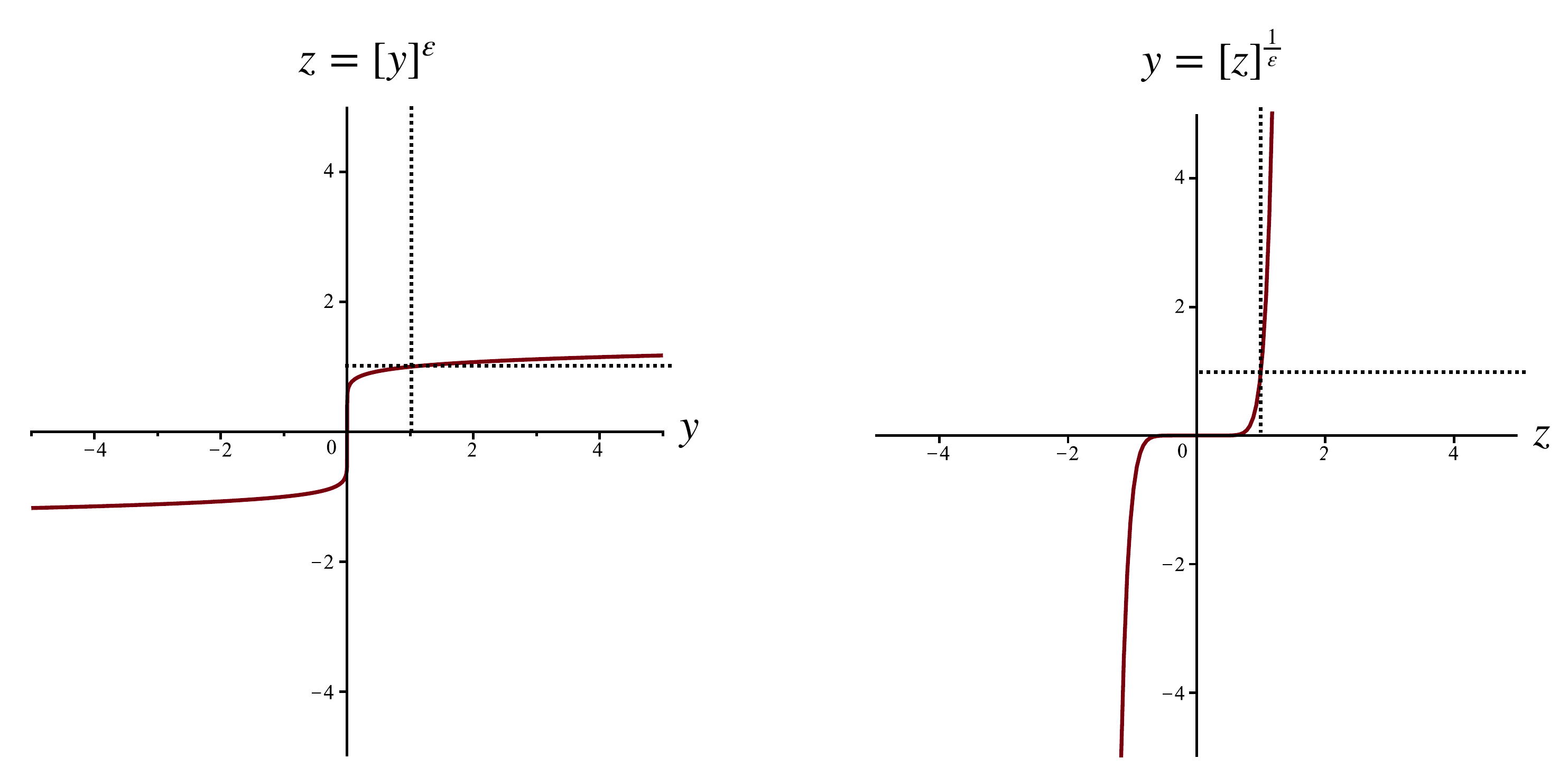}
 \caption{ The function $y \mapsto [y]^\eps$ is strictly increasing from $\Rmat$ into $\Rmat$ ; its reciprocal function is $y = [z]^{\frac{1}{\eps}}$ and we have $\left([y]^{\eps} \right)' = \eps\frac{[y]^{\eps}}{y}$.
When $\eps$ is {infinitely small} this change of variable spreads over $]-1,+1[$ the {halo} of $0$ and concentrates around $+1$ (resp. $-1$) the positive (resp. negative) {\em limited} reals. Here $\eps = 0.1$.} \label{loupe}
 \end{center}
 \end{figure}
 
\noindent If we make the change of variable $z = [y]^{\eps}$ in the system $S_m$ we get:
\beq
\begin{array}{l}
\displaystyle \frac{dz}{dt} =\displaystyle  \eps \frac{[y]^{\eps}}{y}\frac{dy}{dt} \\[8pt]
\displaystyle\frac{dz}{dt} =\displaystyle  \eps \frac{[y]^{\eps}}{y}\left(\frac{1}{\eps}\right) \sqrt{m^2+y^2}(f(x)-y)=\frac{[y]^{\eps}}{y} \sqrt{m^2+y^2}(f(x)-y)\\[8pt]
\displaystyle\frac{dz}{dt} =\displaystyle \frac{[y]^{\eps}}{y} \sqrt{m^2+y^2}(f(x)-y) = z\,\mathrm{sgn}(y) \sqrt{1+\frac{m^2}{y^2}}(f(x)-y)\\[8pt]
\displaystyle\frac{dz}{dt} =\displaystyle |z| \sqrt{1+\frac{m^2}{y^2}}(f(x)-[z]^{\frac{1}{\eps}})
\end{array}
\feq
\textbf{For $ m <1$ we write $m$ in the form:}

$$m = \exp\left(\frac{\rho}{\eps}\right)\quad \quad \rho < 0$$
which gives:
$$\frac{m^2}{y^2} = \exp \left(2\frac{\rho -\ln(|z|)}{\eps}\right)$$so, in the variables $(x,z)$, we study the system :  
\beq \label{eqlocal3}
S_m \quad \quad \quad  \left\{
\begin{array}{lcl}
\displaystyle \frac{dx}{dt}& =&1 \\[6pt]
\displaystyle \frac{dz}{dt} & =&\displaystyle |z|\sqrt{  1+\exp \left(2 \frac{\rho-\ln(|z|)}{\eps} \right)    } \left( f(x)-[z]^{\frac{1}{\eps}}\right)
\end{array} 
\right.
\feq
extended by continuity for $z = 0$.

\bitbul
\item When $z \gnsim 1 $ (resp $z \lnsim -1 $) we have :
	\bito
	\item $|z| > $1 
	\item $\sqrt{\cdots} > $1 
	\item $(f(x) - [z]^{\frac{1}{\eps}}) = -\infty$ (resp $+\infty$) because $[z]^{\frac{1}{\eps}}$ is {infinitely large} positive (resp. negative)
	\fit
and therefore $\displaystyle \frac{dz}{dt} = -\infty$ (resp $+\infty$)

 \item The quantity $\emat^{ \left(2 \frac{\rho-ln(|z|)}{\eps} \right) }$ is {infinitely large} when $\rho -\ln(|z|) \gnsim 0$ and therefore 
 $$ |z|\lnsim \emat^{\rho}\Longrightarrow \frac{dz}{dt} = \mathrm{sgn}(f(x))\cdot(+\infty)$$

\item In the strips\footnote{For the meaning of the double bracket $\eset{}$ see \eqref{ensembleexterne}}
$$\Rmat \times \eset{z\in \Rmat ;\;-1 \lnsim z \lnsim - \emat^{\rho}}\;\mathrm{ and }\; \Rmat \times \eset{z \in \Rmat ; \; \emat^{\rho} \lnsim z \lnsim 1} $$
 we have $\rho+\ln(|z|) \gnsim 0$ and thus $\exp \left(-2 \frac{\rho+\ln(|z|)}{eps} \right) \sim 0$. So in these two strips the trajectories are infinitely close to the trajectories of the vector field: 
 \beq \label{eqlocal4}
\stackrel{\sim}{S_m} \quad \quad \quad \quad \left\{
\begin{array}{lcl}
\displaystyle \frac{dx}{dt}& =&1 \\[6pt]
\displaystyle \frac{dz}{dt} & =& |z]f(x)
\end{array} 
\right.
\feq
 \fit
 From this information it is possible to sketch (see figure \ref{schema2}) the phase portrait of $S_m$ in the variables $(x,z)$.

\begin{figure}[t]
  \begin{center}
 \includegraphics[width=1\textwidth]{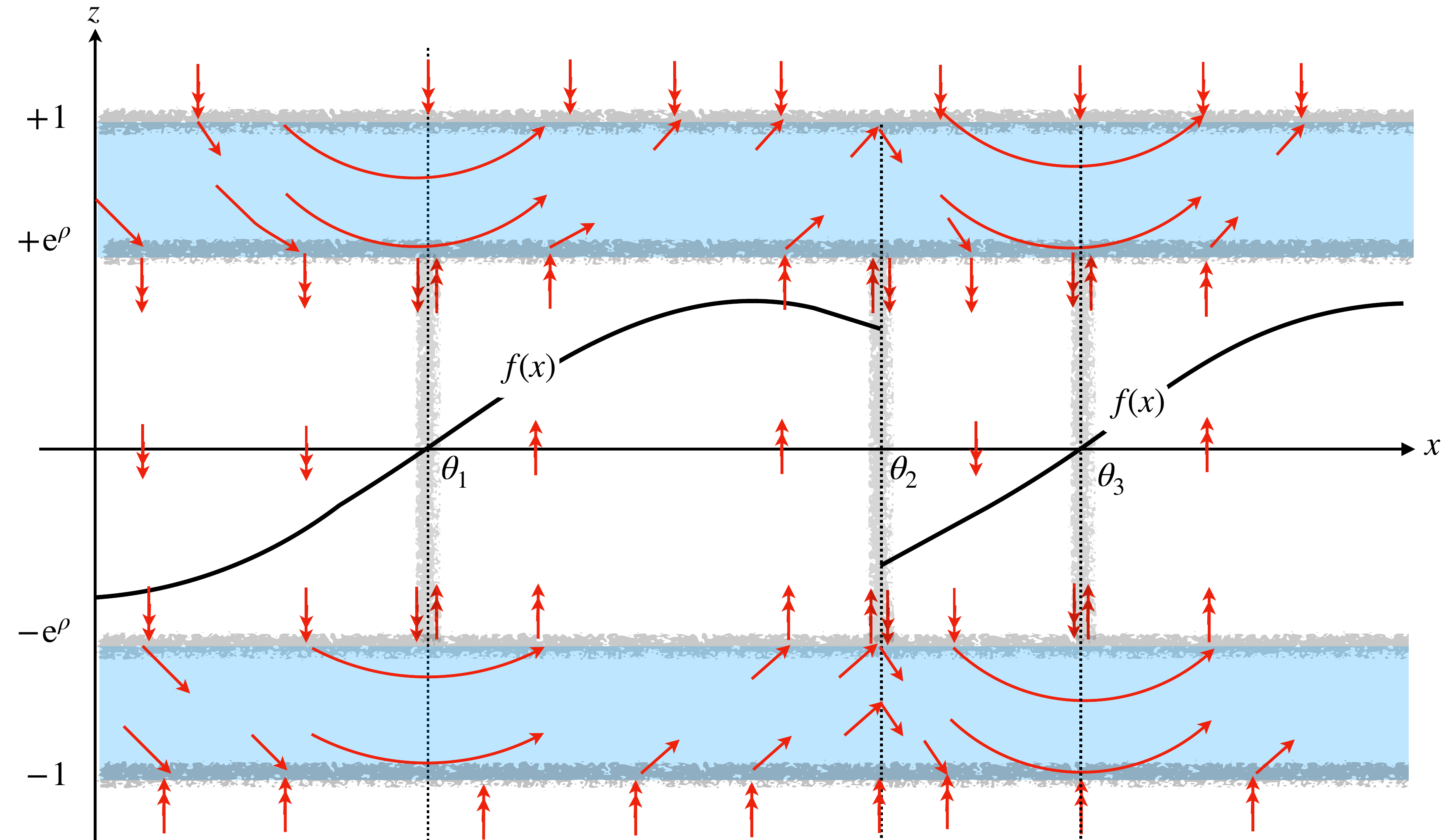}
 \caption{The $S_m$ system under the exponential magnifier. Warning : the graph of $f$ is in the original variables $(x,y)$} \label{schema2}
 \end{center}
 \end{figure}

\subsubsection{Comments on figure \ref{schema2} }

\ben
\item The figure shows a partition in areas delimited by horizontals and verticals.
\item The dotted verticals are the straight lines $x = \theta_n \;n\in J$ where the function $f$ changes sign.
\item The thick shaded horizontals represent the {\em halos} of the lines $z = -1$, $z = -\emat^{\rho}$, $z = \emat^{\rho}$ and $z = 1$ ;
\item The variables are $(x,z)$ but the graph of $f$ is represented in the variables $(x,y)$ ; it is only there to materialize the values of $x$ where the function $f$ changes sign ; on this example there are three changes of sign (for $x = \theta_1,\,\theta_2,\,\theta_3$) ; the change of sign in the middle corresponds to a discontinuity, the two other ones correspond to zeros of $f$. As above, the {halos} of the vertical segments $x = \theta_n$ between $-1$ and $+1$ are symbolized in grey.
\item In the non-grey areas, the second component of the field $S_m$ is {\em infinitely large} positive or negative. This is symbolized by double oriented vertical arrows; above $+1$ they are always pointing down, below $-1$, pointing up or down  in the band $ -\emat^{\rho}< z < \emat^{\rho}$ depending on the sign of $f(x)$.
\item In the blue areas we have symbolized (in red) the variations of the trajectories of the vector field $\stackrel{\sim}{S_m}$ in the vicinity of the zone changes.
\item In the grey areas the trajectory evolves towards the right at unit speed. In particular the vertical grey areas are crossed in an infinitesimal time.
\fen
 \clearpage
 \newpage

\begin{figure}[h]
  \begin{center}
 \includegraphics[width=0.93\textwidth]{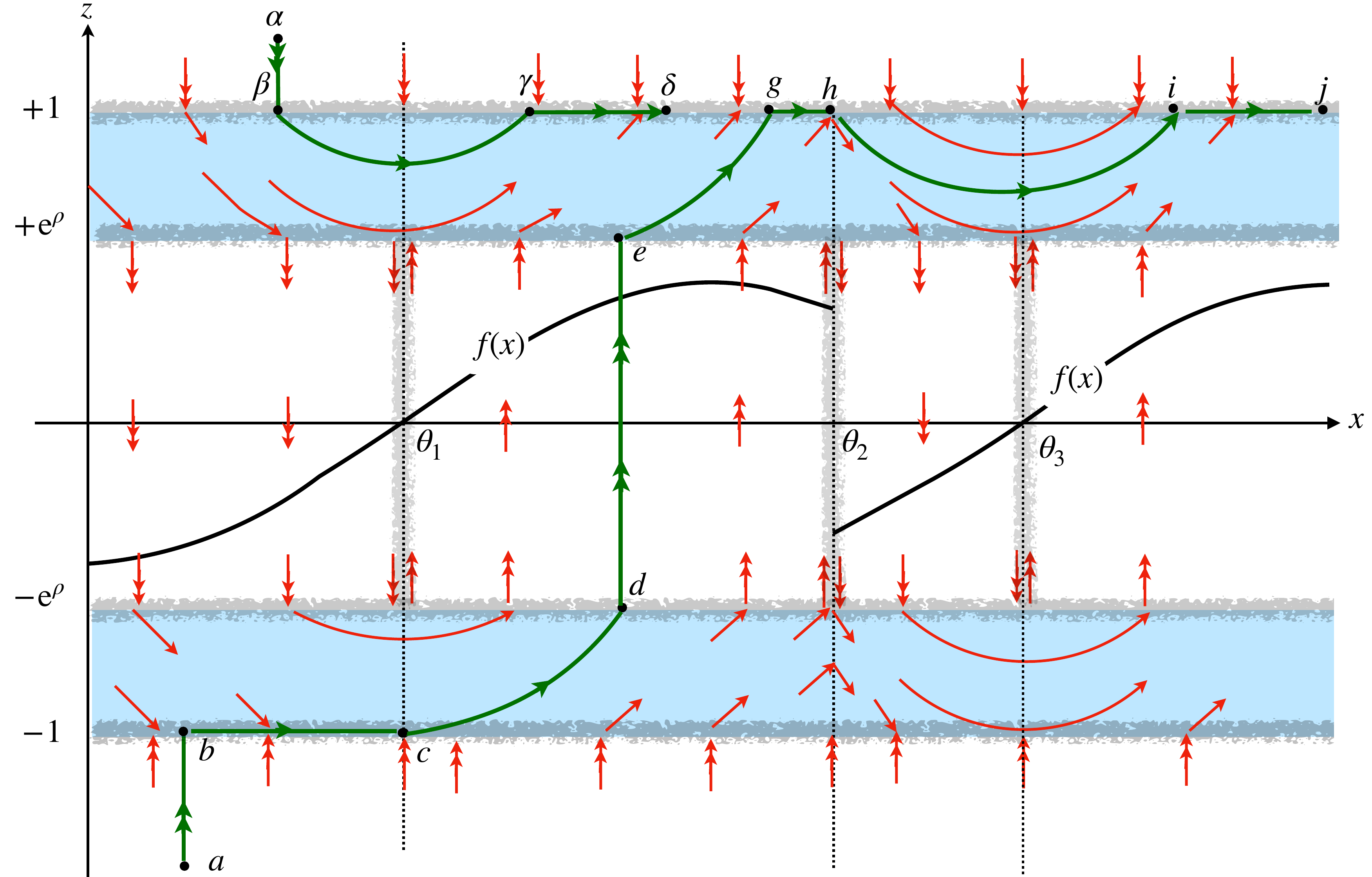}
 \caption{Explanation paragraph \ref{spacexz}.} \label{schema2-2}
 \end{center}
 \end{figure}

 \subsubsection{Comments on figure \ref{schema2-2}: following trajectories in $(x,z)$ variables.}\label{spacexz}
We consider an initial condition $a = (x_0,z_0)$ with $x_0 >0$ and $z_0 \lnsim -1$ not infinitely large, as shown on the figure; the point $(x_0,-1)$ is attractive (i.e. above the field is directed downwards, upwards below); we note $(x(t),z(t))$ the corresponding trajectory. We obviously have $x(t) = x_0+t$.
\paragraph{Trajectory $a,b,c,d,e,g,h,i,j$:}
\ben
\item \texttt{Segment}  $\arc{a,b}$. The field is infinitely large and directed upwards. The corresponding trajectory is almost vertical; there exists a $t _1\sim 0$ such that $z(t_1)\sim-1$ (see appendix \ref{demiLR} for a formal demonstration).

\item \texttt{Segment} $\arc{b,c}$. Since $z = -1$ is attractive up to point $c$ the trajectory {\em hugs} the line $z = -1$ at speed $+1$ up to point $c$ (see appendix \ref{demiLR}) which it reaches at time $t_2 \sim \theta_1-x_0$. 

\item \texttt{Segment} $\arc{c,d}$. There exists an instant $t_3\sim t_2$ such that $(x(t_3),z(t_3))\sim c$ is in the blue zone where $S_m$ is {infinitely close} to $\stackrel{\sim}{S_m}$. So the \texttt{Segment} $\arc{c,d}$ is {infinitely close} to the trajectory $(\tilde{x}(t), \tilde{z}(t))$ of $\stackrel{\sim}{S_m}$ coming from $c$ as long as $\tilde{z}(t) \leq - \emat^{\rho}$ (see appendix \ref{changementdesigne}).
We have:
\beq
\begin{array}{lcl}
\tilde{x}(t )&=& \theta_1 +t\\[6pt]
\tilde{z}(t) &=& \displaystyle -\emat^{- \int _{\theta_1} ^t f(s) ds} 
\end{array}
\feq
and so the value $-\emat^{\rho}$ is reached at time $\tau_4$ such that :
\beq 
\displaystyle -\emat^{- \int _{\theta_1} ^{\tau_4}f(s) ds} = -\emat^{\rho}
\feq 
or :
\beq \label{tau1}
\displaystyle \int _{\theta_1} ^{\tau_4}f(s) ds = -\rho
\feq 
Thus there exists a $t_4 \sim \tau_4$ such that $ z(t_4) \sim -\emat^{\rho} $

\item \texttt{Segment} $\arc{d,e}$. From $t_4$ the trajectory is quasi-vertical ascending. There exists $t_5 \sim t_4$ such that $z(t_5) \sim \emat^{\rho} $. At point $e$ the field crosses the line $z = +1$. 

\item \texttt{Segment} $\arc{e,g}$. From $t_5$ we are in the same situation as at time $t_3$ (see 3. above). The trajectory of $\stackrel{\sim}{S_m}$ thus reaches the value $1$ at time $\tau_6$ such that :
\beq 
\displaystyle \emat^{\rho} \emat^{ \int _{\tau_5} ^{\tau_6}f(s) ds} = 1
\feq 
or :
\beq \label{tau2}
\displaystyle \int _{\tau_5} ^{\tau_6}f(s) ds = -\rho
\feq 

\item By joining \eqref{tau1} and \eqref{tau2} it comes that $t_6 \sim \tau_6$ where $\tau_6$ is defined as the first instant such that :
\beq \label{condition1}
\displaystyle \int_{\theta_1} ^{\tau_6}f(s) ds = - 2\rho
\feq 

\item \texttt{Segment} $\arc{g,h}$. For the same reasons as 2. above the segment $\arc{f,g}$ is in the {halo} of $z = +1$ as long as $t \lnsim \theta_2$.
\item \texttt{Segment} $\arc{h,i}$. For the same reasons as 3. above the trajectory of $S_m$ remains in the blue zone until the time $t_7 \sim \tau_7$ defined by :
\beq \label{condition2}
\displaystyle \int_{\theta_2} ^{\tau_7}f(s) ds = 0
\feq 
Unlike the segment $\arc{c,d}$ which crossed the blue zone between $z = - 1$ and $z = -\emat^ {\rho} $ here the segment $\arc{h,i}$, originating from a point such as $z = 1$, goes back up to a point such as $z = 1$ without crossing the blue zone.
\item \texttt{Segment} $\arc{i,j}$. As for $\arc{b,c}$ we stay in the {\em halo} of $z = 1$.
\fen
\paragraph{Trajectory $\alpha,\beta,\gamma,\delta$:}
While the previous trajectory crosses the band $-1< z< +1$ this trajectory enters the blue area but does not cross the band nor the axis $y = 0$.

\begin{remarque}
There is no qualitative difference between the point $g$ which corresponds to a discontinuity of $f$ and the point $c$ which corresponds to a zero of $f$.
\end{remarque}

 \subsection{Following a trajectory : \\Simultaneous use of the variables $x,y$ and $x,z$.}
\begin{figure}[h]
  \begin{center}
 \includegraphics[width=1\textwidth]{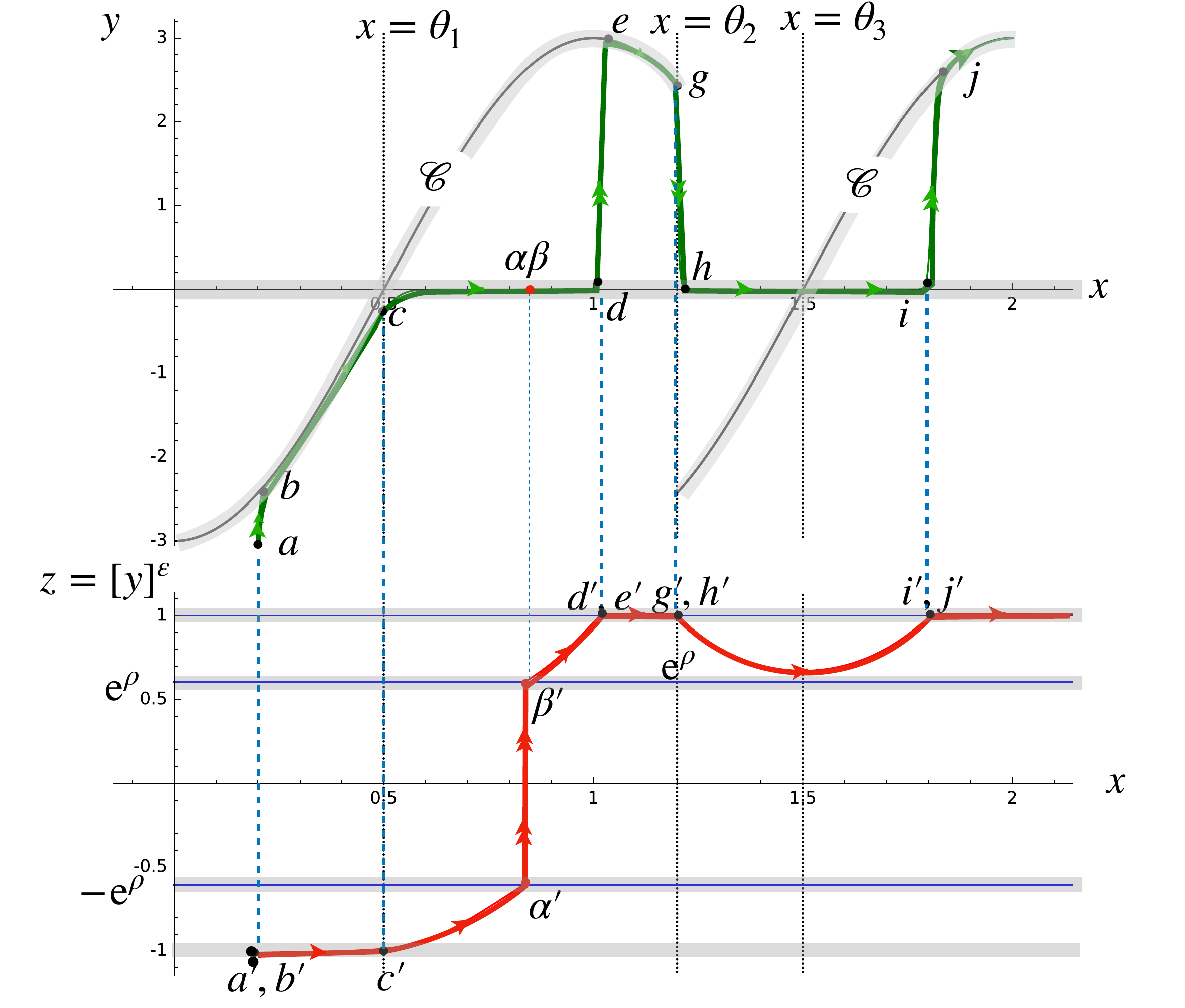}
 \caption{Above the system $S_m$ in the variables $x,y$ and below the enlargement of the halo of the axis $y=0$ under the exponential magnifying glass.} \label{loupesimu1}
 \end{center}
 \end{figure}
 
 On the figure \ref{loupesimu1} are superimposed the spaces of the trajectories of $S_m$, (with $m= \emat^{\frac{\rho}{\eps}}$) in the variables $(x,y)$ and the variables $(x,z)$, the values of $x$ corresponding to each other. 
 \bit
 \item In the variables $(x,y)$ we have represented the graph (curve $\mathcal{C}$) of the function $f$ used in the previous section, which cancels twice at $x =\theta_1$ and $x =\theta_3$ and has a discontinuity where it changes sign at $x =\theta_2$. The halos of $\mathcal{C}$ and the axis $y = 0$ are symbolized in grey.
  \item In the variables $(x,z)$ the halos of the horizontals $z = ±1$ and $z = ± \emat^ {\rho}$ are symbolized in gray.
 \fit 
 We will follow the trajectory from point $a$.
 \bitbul
 \item \texttt{Segment $\arc{a,b}$}. The slow curve is attractive. We follow a vertical segment and then enter the halo of $\mathcal{C}$ at a point $b$ with abscissa infinitely close to that of $a$. In the variables $(x,z)$ the two corresponding points, $a'$ and $b'$ are in the halo of $z = - 1$ and infinitely close to each other.

 \item\texttt{Segment $\arc{b,c}$}. The trajectory follows the attractive slow curve as long as $f(x) \lnsim 0$ until the point $c = (x_1,y_1)$ where it enters the halo of $y = 0$. In space $(x,z)$ the segment $\arc{b',c'}$ remains in the halo of $z = -1$ until the point $c'$ infinitely close to $(\theta_1,-1)$.
\fit 
In the space $(x,y)$, from the point $c$ we have two possibilities: to stay in the halo of $\mathcal{C}$ or that of $y = 0$; the answer is given in the space $(x,z) $. We now follow the trajectory in $(x,z)$ space from the point $c'$ to the point $d'$.
\bitbul
 \item\texttt{Segment $\arc{c',\alpha'}$}.  We now follow the trajectory in $(x,z)$ space from the point $c'$ to the entry into the halo of the line $z = -\emat^{\rho}$ at the point $\alpha'$. In $(x,y)$ space, since $|z| \lnsim 1$, the segment $\arc{c,\alpha}$ remains in the halo of $y = 0$.

 \item\texttt{Segment $\arc{\alpha',\beta'}$}. In the space $(x,z)$ the segment $\arc{\alpha',\beta'}$ makes jump in an infinitely small time from $z =-\emat^{\rho}$ to $z =+\emat^{\rho}$ (paragraph \ref{spacexz}, point 3. ) In space $(x,y)$ the segment $\arc{\alpha,\beta}$ is infinitely short and crosses the axis $y = 0$.
 
  \item \texttt{Segment $\arc{\beta',d'}$}. We follow the trajectory from the point $\beta'$ until it meets the halo of the line $z = +1$ at the point $d'$. To enter the halo of $z = 1$ is to leave the halo of $y = 0$, therefore, in the space $(x,y)$ the corresponding point $d = (x_2,y_2) $ is the point of exit of the halo of $y = 0$. The value of $x_2$ is given by :
  $$\int_{\theta_1}^{x_2}f(s)ds =- 2\rho$$ 
  
  \fit
    Since we are out of the halo of $y = 0$ we can go back to the space of variables $(x,y)$.
  \bitbul
  
 \item\texttt{Segment $\arc{d,e,g,h}$}. We jump to $e$ where we enter the halo of the $\mathcal{C}$ curve which we follow until the point $g$ where we jump to the halo of $y =0$ where we enter the point $h$. The corresponding segment in the space $(x,z)$, $\arc{d',e,',g',h'}$ runs along $z = 1$.
 \fit
 We return to the space $(x,z)$.
 \bitbul
 
  \item \texttt{Segment $\arc{h',i'}$}. In $(y,z)$ space we follow the trajectory that does not cross the $ \emat^{\rho} < z < 1$ band but joins the $z =1$ halo again at the  point $i'$. In $(x,y)$ space, the corresponding point $h =(x_3,y_3)$ is the exit point from the halo of  $y = 0$ . The value of $x_3$ is given by:
  $$\int_{\theta_1}^{x_3}f(s)ds = 0$$
 
 \fit
 We return to the space $(x,y)$.
 \bitbul
  \item \texttt{Segment $\arc{i,j}$}. In $(x,y)$ space we jump from $i$ to point $j$ in the halo of $\mathcal{C}$.
 \fit

 \subsubsection{Initial condition in the halo of the $y=0$ axis}\label{ci0}
  \begin{figure}[h]
  \begin{center}
 \includegraphics[width=1\textwidth]{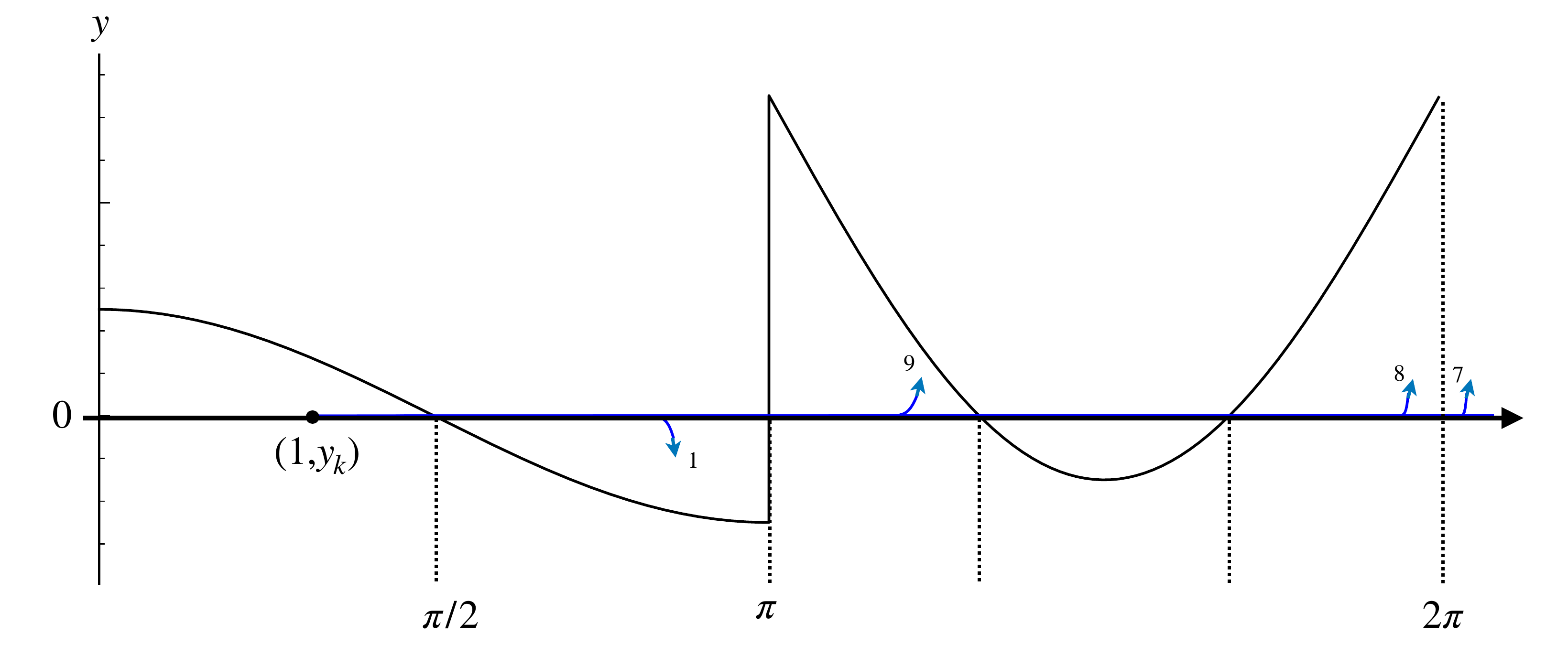}
 \caption{Simulation of 9 trajectories from 9 infinitely close initial conditions. Explanations paragraph \ref{ci0} } \label{cisim01}
 \end{center}
 \end{figure}
  \begin{figure}[h]
  \begin{center}
 \includegraphics[width=1\textwidth]{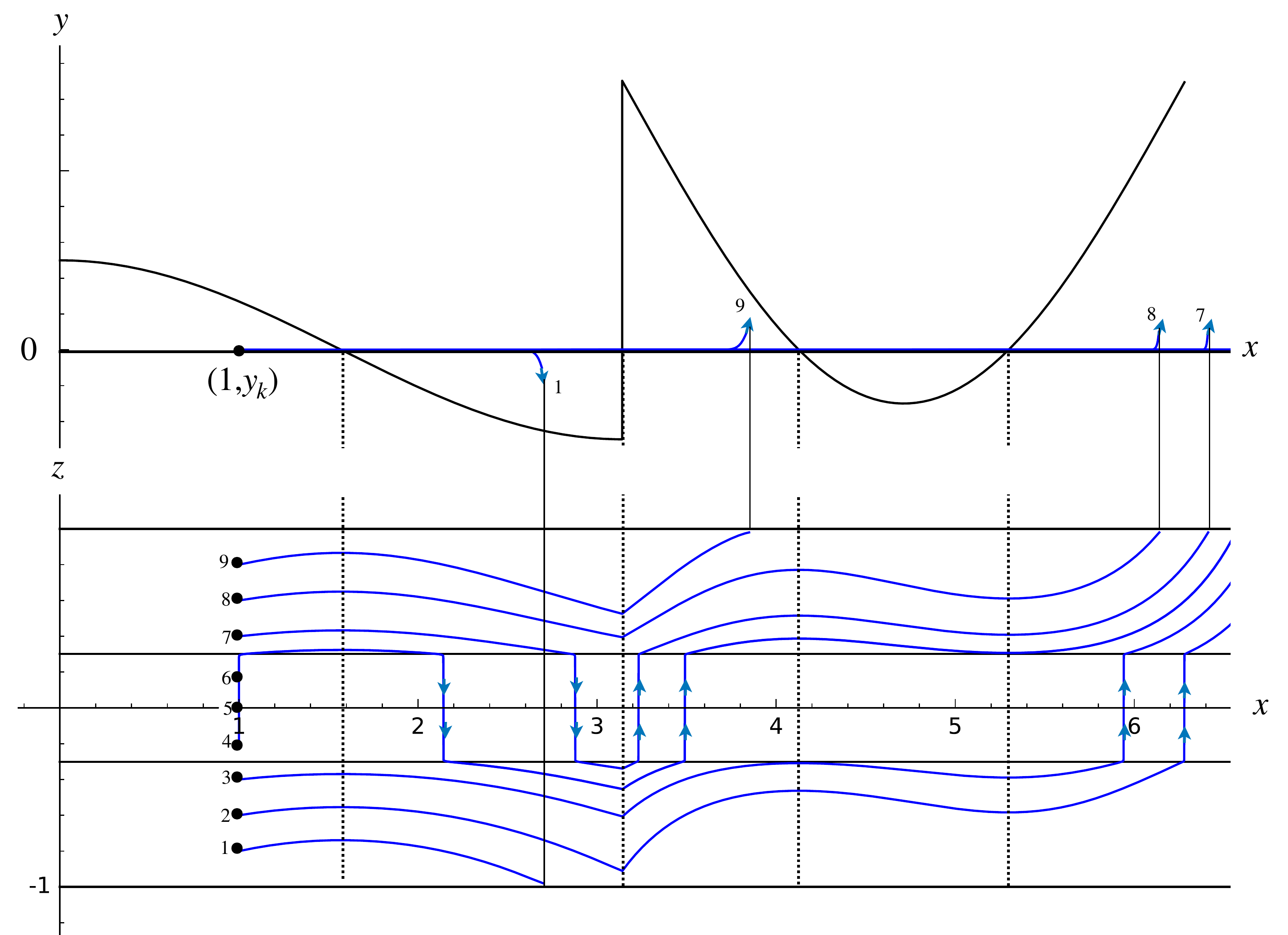}
 \caption{Magnification of the simulation of the trajectories of simulation \ref{cisim01}. Explanations paragraph \ref{ci0} } \label{cisim02}
 \end{center}
 \end{figure}
 Let$(x_0, y_0)$ be an initial condition  such that $y_0 \sim 0$ (condition excluded by the approximation theorem). What can we say about the solution from $(x_0, y_0))$ without any other information about $y_0$ ? 
 
 On the figure \ref{cisim01} we observe the result of the simulation of $S_m$ with the function $f$ given by :
 \beq \label{fcisim0}
\begin{array}{lcl}
\displaystyle x \in [0,\pi[ &\Longrightarrow &f(x) = 0.5\cos(x) \\[8pt]
\displaystyle x \in [\pi,2 \pi[& \Longrightarrow &f(x) = 1.5 + 1.8\sin(x)
\end{array}
\feq
 whose graph is plotted in black, for the parameter values :
 $$ \eps = 0.01 \quad \quad m = \emat^ {\rho/\eps} \quad \rho = -1.2$$
 and from the $9$ initial conditions:
 $$ (x_0 = 1, y_k)\quad y_k = (k-5)\times 0.2 \times m \quad k = 1.2,\cdots, 9$$
 on a duration of 5.5 units. The initial conditions are very close; with the parameter values we have $|y_k|< 2\times 10^{-10}$, the 9 initial conditions are thus indistinguishable on the figure.
  The simulation shows that on the $9$ trajectories, only $4$ (corresponding to $k = 1,7,8,9$) have left the halo of $y = 0$ at different points before the time $t =5.5$ .
 
 The use of the exponential magnifying glass (in the simulation instead of displaying $y$ we display $z = \emat^{\eps\log(y)}$ if $y >0$ and $z = -\emat^{\eps\log(-y)}$ if $y <0$) separates the initial conditions and the associated trajectories. This is what is shown on  figure \ref{cisim02} where we see that the exit point of the halo of the $y = 0$ axis depends in an essential way on the (infinitely small) distance of $y_k$ from $0$. This is the reason why, in the approximation theorem, the points of the halo on the $y = 0$ axis are excluded as possible initial conditions. On the other hand, when a trajectory enters the halo of the $y = 0$ axis for a value close to $x$, it can only do so by passing through the values $z \sim ± 1$ which removes any uncertainty on the value $S_{\rho}(x)$ of the exit point.

 \subsection{Foliation of the trajectories when $f$ changes  sign}
  \begin{figure}[h]
  \begin{center}
 \includegraphics[width=1\textwidth]{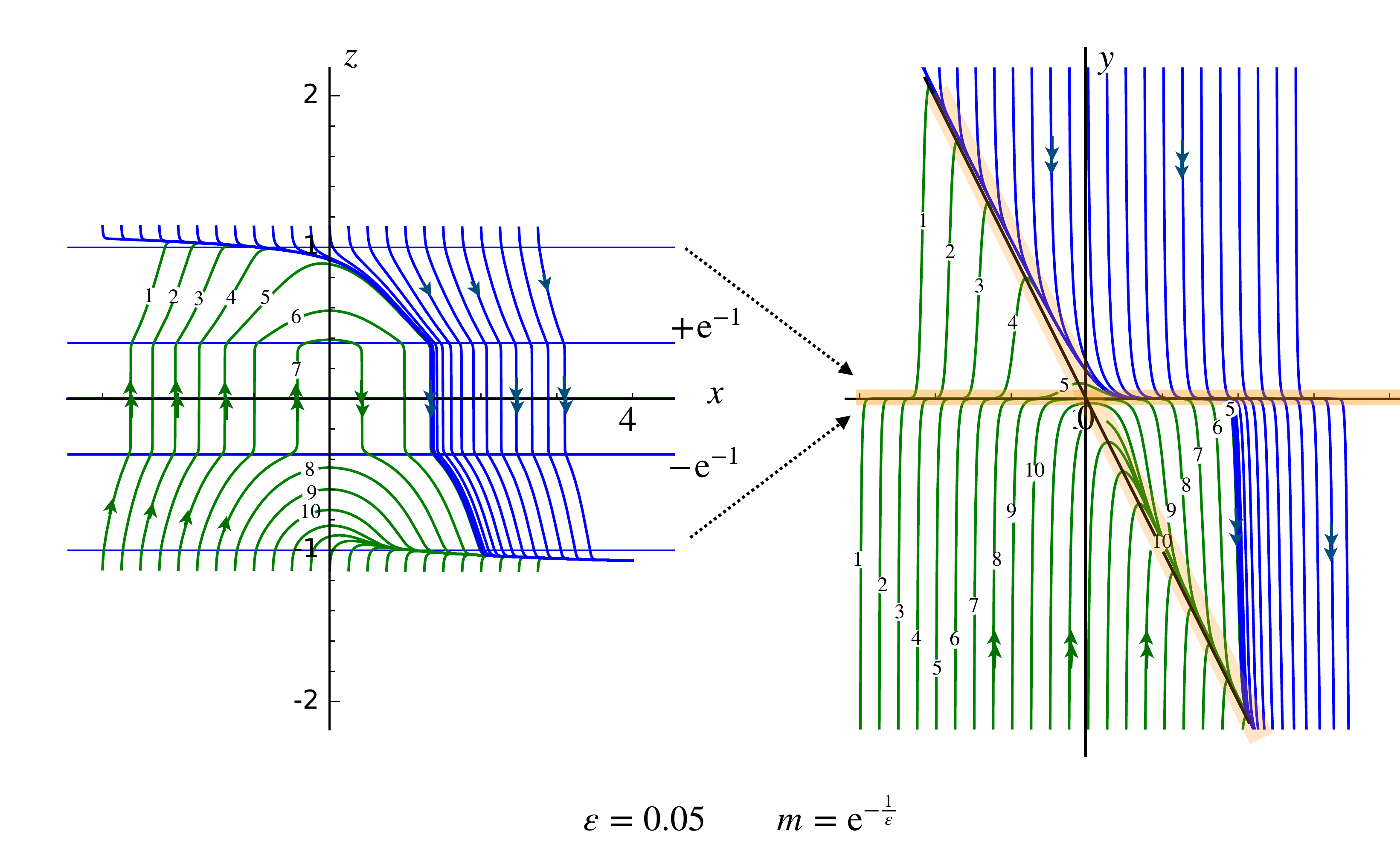}
 \caption{Simulation of $S_m$ with $f(x) = -x$} \label{retard4}
 \end{center}
 \end{figure}
  \begin{figure}[h]
  \begin{center}
 \includegraphics[width=1\textwidth]{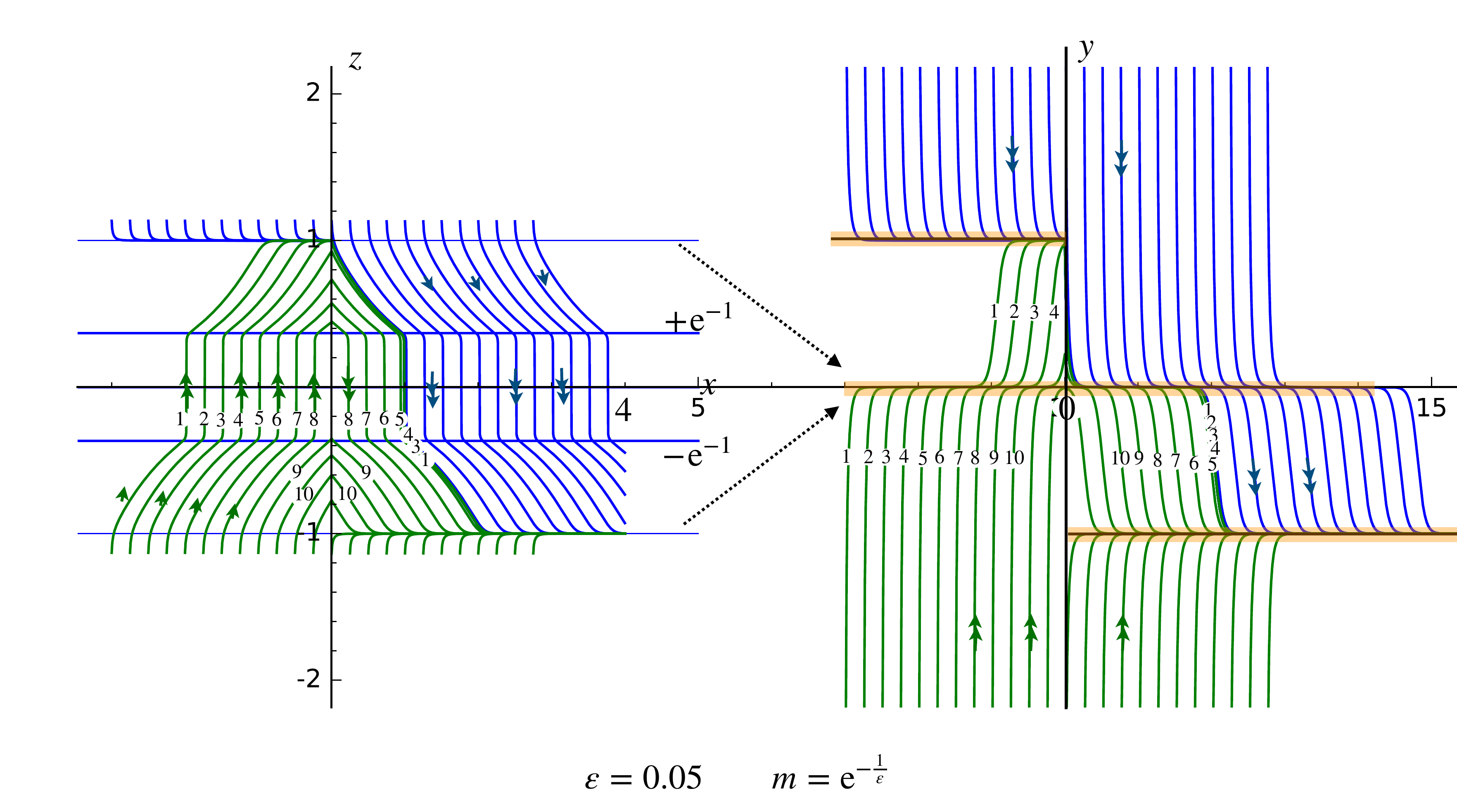}
 \caption{Simulation of $S_m$ with $f(x) = ±1$} \label{retard4bis}
 \end{center}
 \end{figure}
 On the figures \ref{retard4} and \ref{retard4bis} we have simulated the trajectories of the system $S_m$ around $(0,0)$ in the case where $f(x) = -x$ and where $f(x) = 1$ if $x$ is negative, $-1$ otherwise, to illustrate the passage of a change of sign of $f$ when $f$ cancels or is discontinuous. On the right we have the trajectories in the original variables where the halos of $y = 0$ and $y = f(x)$ are symbolized by the yellow bands and on the left under the exponential magnifier.
 \clearpage
\newpage

 \section{Proof of a conjecture of G. Katriel}
             
\subsection{The problem }           
            
             We consider the system :
             
 \beq \label{Kat1}
\begin{array}{lcl}
\displaystyle \frac{dx_1}{dt} &=&r_1(\nu t)x_1+mu(x_2-x_1)\\[6pt]
\displaystyle \frac{dx_2}{dt} &=& r_2(\nu t)x_1+\mu(x_1-x_2)
 \end{array} 
\feq
where the functions $r_i$ are standard, periodic of period $2\pi$ and piecewise differentiable in the sense of the section \ref{equationsdebase}. We are interested in the behavior of the quantity :
\beq \label{Delta}
\displaystyle \Delta(\nu,\mu) = \lim_{t \rightarrow \infty} \frac{\ln(x_1(t)x_2(t))}{t} 
\feq
when :
\beq \label{HK}
\displaystyle \frac{1}{2\pi} \int_0^{2\pi} r_i( t)dt < 0 \quad \chi = \frac{1}{2\pi} \int_0^{2\pi} \max_{i = 1,2} r_i(t)dt >0
\feq
The system \eqref{Kat1} represents the growth of two populations of size $x_1$ and $x_2$, at two sites $1$ and $2$, supporting to different environments, $r_i(t)$, in the presence of migration between the two sites. The assumptions on $r_i$ have the consequence that, in the absence of migration ($\mu = 0$), at each site the populations disappear. Let's call {\em inflation threshold} and note $\mu^*(\nu)$ :
\beq \label{threshold }
\displaystyle \mu^*(\nu) = \inf_{\mu \geq 0} \mu : \Delta(\nu,\mu) > 0 \}
\feq

In  \cite{KAT21} it is conjectured that :
\begin{proposition}\label{conjecturekat} When $\nu$ is infinitely small (the period of variation of the environment is very large) the threshold of inflation is of the order of $\exp(-\frac{1}{\nu})$ (exponentially small compared to $\nu$).
\end{proposition} 
I demonstrate this conjecture in the next paragraph. I refer to \cite{BLSS21, HOLTPNAS20, KLA08} for more information on the inflation phenomenon highlighted in \cite{HOLT02}.

\subsection{Reduction to a system of the type $S_m$} 
\paragraph{Change to the variables $U = \ln(x_1x_2)$ and $V = \ln(x_1/x_2)$} $\,$ \\
Let's put:
$$ \xi_1 = \ln(x_1)\quad \xi_2 = \ln(x_2)$$
In these variables the system \eqref{Kat1} becomes :
  \beq \label{Kat2}
\begin{array}{lcl}
\displaystyle \frac{d\xi_1}{dt} &=&r_1(\nu t)+\mu\left(\emat^{\xi_2-\xi_1}-1\right)\\[6pt]
\displaystyle \frac{d\xi_2}{dt} &=& r_2(\nu t) )+\mu\left(\emat^{\xi_1-\xi_2}-1\right) 
 \end{array} 
\feq
and now we pose :
$$ U = \xi_1+\xi_2\quad \quad V = \xi_1 - \xi_2$$
In these variables the system \eqref{Kat2} becomes :
 \beq \label{Kat3} 
\begin{array}{lcl}
\displaystyle \frac{dU}{dt} &=&\displaystyle r_1(\nu t)+r_2(\nu t) + 2\mu (\ch(V)-1) \\[8pt]
\displaystyle \frac{dV}{dt} &=&(r_1(\nu t)-r_2(\nu (t))-2\mu\,\sh(V)
 \end{array} 
\feq
We see that the variable $V$ is decoupled from $U$. We have:
$$\Delta(\nu,\mu) = \lim_{t \rightarrow \infty} \frac{U(t)}{t} $$
Consider the equation :
\beq \label{V1}
  \frac{dV}{dt} =(r_1(\nu t)-r_2(\nu (t))-2\mu\,\sh(V)
\quad \quad 
\feq
It is a periodic non-autonomous equation.
\begin{proposition}
When $\mu > 0$ the equation \eqref{V1} has a unique periodic solution (of period $\frac{2\pi}{\nu}$) globally asymptotically stable $V_{\nu,\mu}(t)$.
\end{proposition}
\textbf{Proof} For $M$ sufficiently large we have for $|V| = M$, $ \displaystyle \mathrm{sgn}\left(\frac{dV}{dt} \right) =-\mathrm{sgn} (V)$ which means that $[-M,+M]$ is invariant, the application $V_0 \mapsto V(2\pi/\nu, 0,V_0)$ has at least one fixed point.
For such a fixed point, the equation of variations along the associated periodic solution $V_{\nu,\mu}(t)$ is :
$$\frac{d\delta V(t)}{dt} = -2\mu\cosh(V_{\nu,\mu} (t))\delta V(t)$$
So the derivative of $V_0 \mapsto V(2\pi/\nu, 0,V_0)$ is positive strictly smaller than $1$ which proves the uniqueness and attractiveness of the fixed point.
$\Box$\\
It follows that:
$$\Delta(\nu,\mu) = \lim_{t \rightarrow +\infty} \frac{U(t)}{t} = \frac{\nu}{2\pi} \int_0^{2\pi/\nu}\displaystyle r_1(\nu s)+r_2(\nu s) + 2\mu \left(\ch(V_{\nu,\mu}(s))-1 \right) ds $$

\paragraph{ Reduction to a fixed period.} The system is reduced to the fixed  period $2\pi$ by putting  $\overline{V} (t) = V(t/\nu)$ which gives :
\beq \label{V2}
  \frac{d\overline{V}}{dt} =\frac{1}{\nu} \big(r_1(t)-r_2(t)-2\mu\,\sh(\overline{V})\big)
\quad \quad 
\feq
and : 
$$\Delta(\nu,\mu) = \frac{1}{2\pi} \int_0^{2\pi} \displaystyle r_1( s)+r_2(s) + 2\mu \left(\ch(\overline{V}_{\nu,\mu}(s))-1 \right) ds$$
where $\overline{V}_{\nu,\mu}$ is the periodic solution of \eqref{V2}.

\paragraph{Linearization of the hyperbolic sine.} We fix $\nu = \eps \sim 0$ (to recall that $\nu$ is small en to fit with the previous notations) once and for all and $\mu$ is a parameter. We define:
$$ W = 2\mu\sinh(\overline{V})$$
which gives :
\beq \label{W1}
\frac{dW}{dt} = \frac{2\mu}{\eps} \cosh(\overline{V}) \big(r_1( t)-r_2(t)-W)\big)
\feq
\beq \label{W2}
\frac{dW}{dt} = \frac{2\mu}{\eps} \sqrt{1+ \sinh^2(\overline{V})} \sqrt{1+ \sinh^2(\overline{V})} \big(r_1( t)-r_2(t)-W)\big)
\feq

\beq \label{W3}
\frac{dW}{dt} = \frac{2\mu}{\eps} \sqrt{ 1+( \sinh(\sinh^{-1}(W/\mu))^2 } \big(r_1( t)-r_2(t)-W\big)
\feq

\beq \label{W4}
\frac{dW}{dt} = \frac{1}{\eps} \sqrt{ 4\mu^2+W^2 } \big(r_1( t)-r_2(t)-W\big)
\feq
and : 
$$\Delta(\eps,\mu) = \frac{1}{2\pi} \int_0^{2\pi}\displaystyle r_1( s)+r_2(s) + 2\mu\left( \sqrt{1+(W_{\mu}/2\mu)^2}-1 \right) ds$$
where $W_{\mu}$ is the periodic solution of \eqref{W4}. If we put $m = 2\mu$ and $W = y$ we finally get :
\beq \label{W5}
\frac{dy}{dt} = \frac{1}{\eps} \sqrt{ m^2+y^2 } \big(r_1( t)-r_2(t)-y\big)
\feq
which, if we add $\frac{dx}{dt} = 1$, is the system $S_m$ with $f = r_1-r_2$. 
So we have:
$$\Delta(\eps,\mu ) = \Delta(\eps,m) = \frac{1}{2\pi} \int_0^{2\pi} r_1( s)+r_2(s) + \left( \sqrt{m^2+y_m^2(s)}-m \right) ds $$
where $y_m$ is the periodic solution of \eqref{W5}.\\\\
We see that if $m = \sim 0$ we have :
$$\Delta(\eps,m) \sim \frac{1}{2\pi} \int_0^{2\pi}\displaystyle r_1( s)+r_2(s) + |y_m(s)| ds $$
The condition $\chi > 0$ of \eqref{HK} can be read:
$$2 \chi = \frac{1}{2\pi} \int_0^{2\pi} 2 \max_{i = 1,2} r_i( s)ds = \frac{1}{2\pi} \int_0^{2\pi}r_1(s)+r_2(s) + |r_1(s)-r_2(s)| ds>0 $$
So proving the proposition \ref{conjecturekat} is the same as showing that, if we put $m = \emat^{\rho/\eps} $ the quantity:
$$ \frac{1}{2\pi} \int_0^{2\pi} |r_1(s)-r_2(s)| - y_m(s) ds$$
can be made smaller than the standard quantity $\chi > 0$ for non-infinitely small values of $\rho$.

\paragraph{Analysis of the periodic solution $y_m$.} To fix the ideas we do this for a particular system; the general analysis obviously proceeds in the same way.

\begin{figure}[h]
  \begin{center}
 \includegraphics[width=1.0\textwidth]{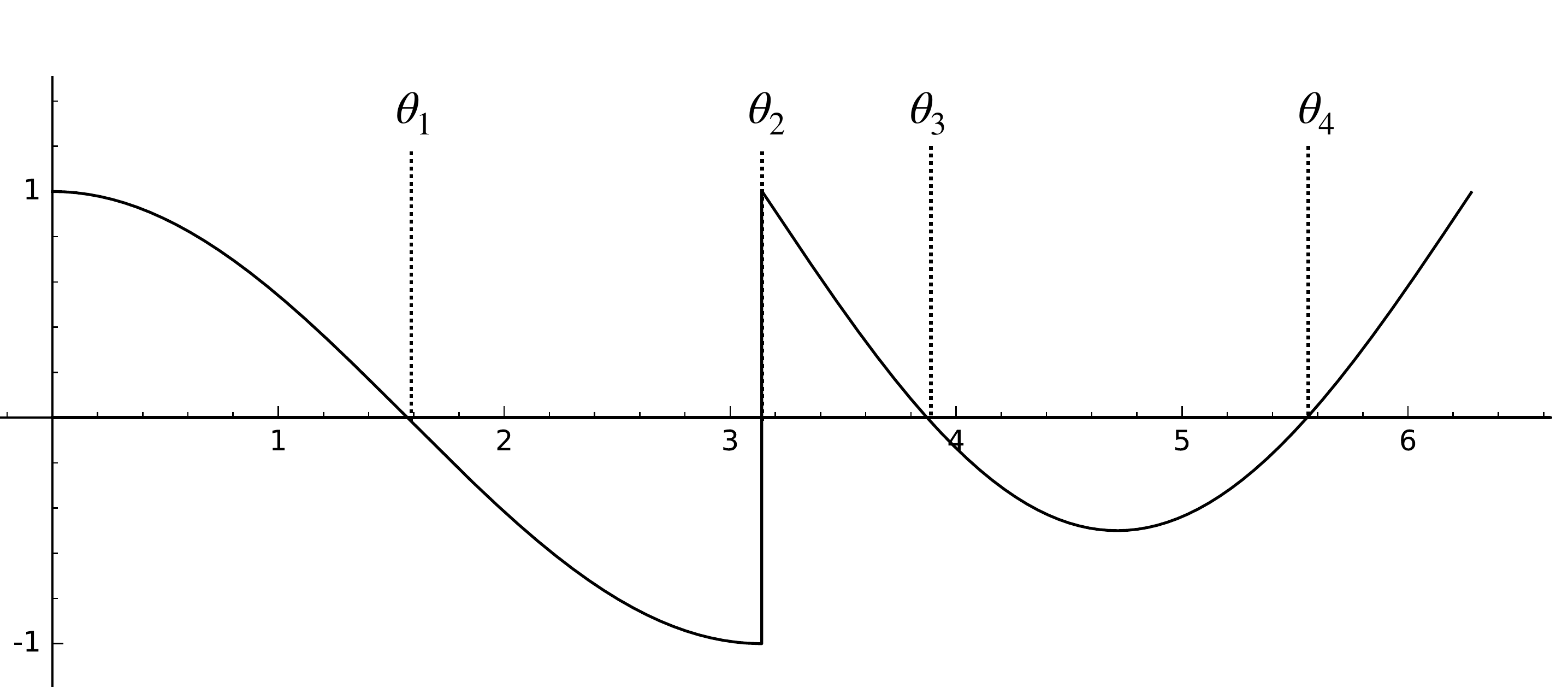}
 \caption{Graph of $r_1-r_2$} \label{grapher1r2}
 \end{center}
 \end{figure}

We consider the system :
\beq \label{per1}
S_m \quad \quad \left \{
\begin{array}{lcl}
\displaystyle \frac{dx}{dt} &=&1\\[6pt]
\displaystyle \frac{dy}{dt}& =&\displaystyle \frac{1}{\eps} \sqrt{ m^2+y^2 } \big(r_1( t)-r_2(t)-y\big)
\end{array}
\right.
\feq
with $r_1-r_2$ of period $2\pi$ defined by : :
\beq \label{r1mr2}
\left[  
\begin{array}{lcl}
\displaystyle x \in [0,\pi[ &\Longrightarrow &r_1(x)-r_2(x) = \cos(x)\\[4pt]
\displaystyle x \in [\pi,2 \pi[& \Longrightarrow &r_1(x)-r_2(x) = 1 + 1.5 \sin(x)
\end{array}
\right.
\feq
This function is continuous at all points, except $x = \pi$, it changes sign for the four values $\theta_1, \theta_2, \theta_3, \theta_4$ (see figure \ref{grapher1r2}).

For $m >0$, $S_m$ has a unique globally asymptotically stable periodic solution $y_m$ which we now analyze for exponentially small values of $m$ of the form $m = \emat^{\rho/\eps}\,\;\rho \lnsim 0$ .  The phase space of \eqref{per1} is $S^1\times \Rmat$ where $S^1 = \Rmat (\mathrm{mod} 2\pi)$. 

To the system $S_m$ we associate the constrained  system  (see section \ref{contraint}):
\beq \label{Cper1}
S_m^0 \quad \left \{
\begin{array}{lcl}
\displaystyle \frac{dx}{dt} &=&1\\[8pt]
\displaystyle 0 & =&\displaystyle \sqrt{ m^2+y^2 } \big(r_1( t)-r_2(t)-y\big)
\end{array}
\right.
\feq

\begin{lemme}
There exists $\rho^* \lnsim 0$ such that for all $\rho \geq \rho^*$ we have : $$S(\theta_i) < \theta_{i+1}\quad (\mathrm{mod},2\pi)$$
\end{lemme}
\textbf{Proof} This is an immediate consequence of the definition of $S(x)$ (see proposition \ref{rderho})\\\\
\noindent From now on we assume that :
$$ \rho^*\leq \rho \lnsim 0$$
  \begin{figure}[h]
  \begin{center}
 \includegraphics[width=1.0\textwidth]{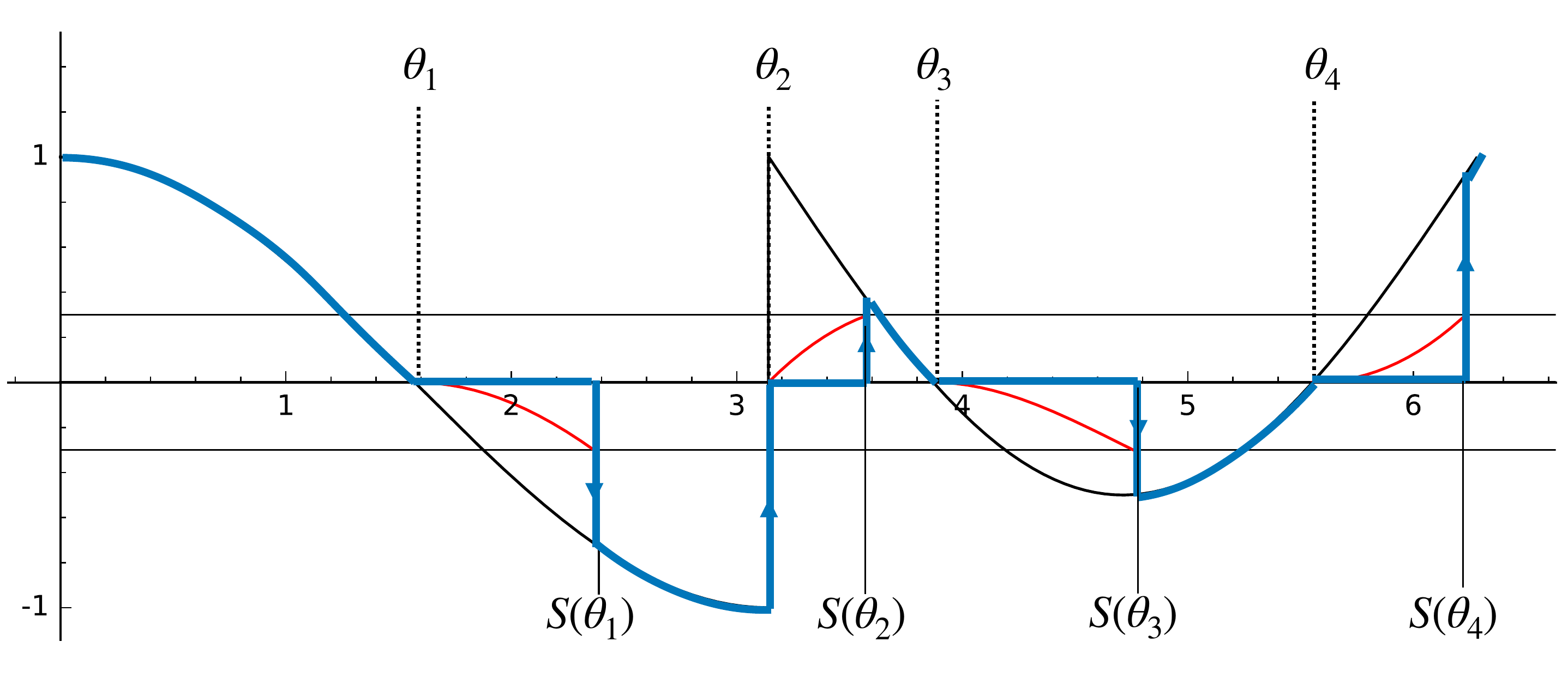}
 \caption{C-periodic trajectory} \label{grapher1r2bis}
 \end{center}
 \end{figure}
  \begin{figure}[h]
  \begin{center}
 \includegraphics[width=1.0\textwidth]{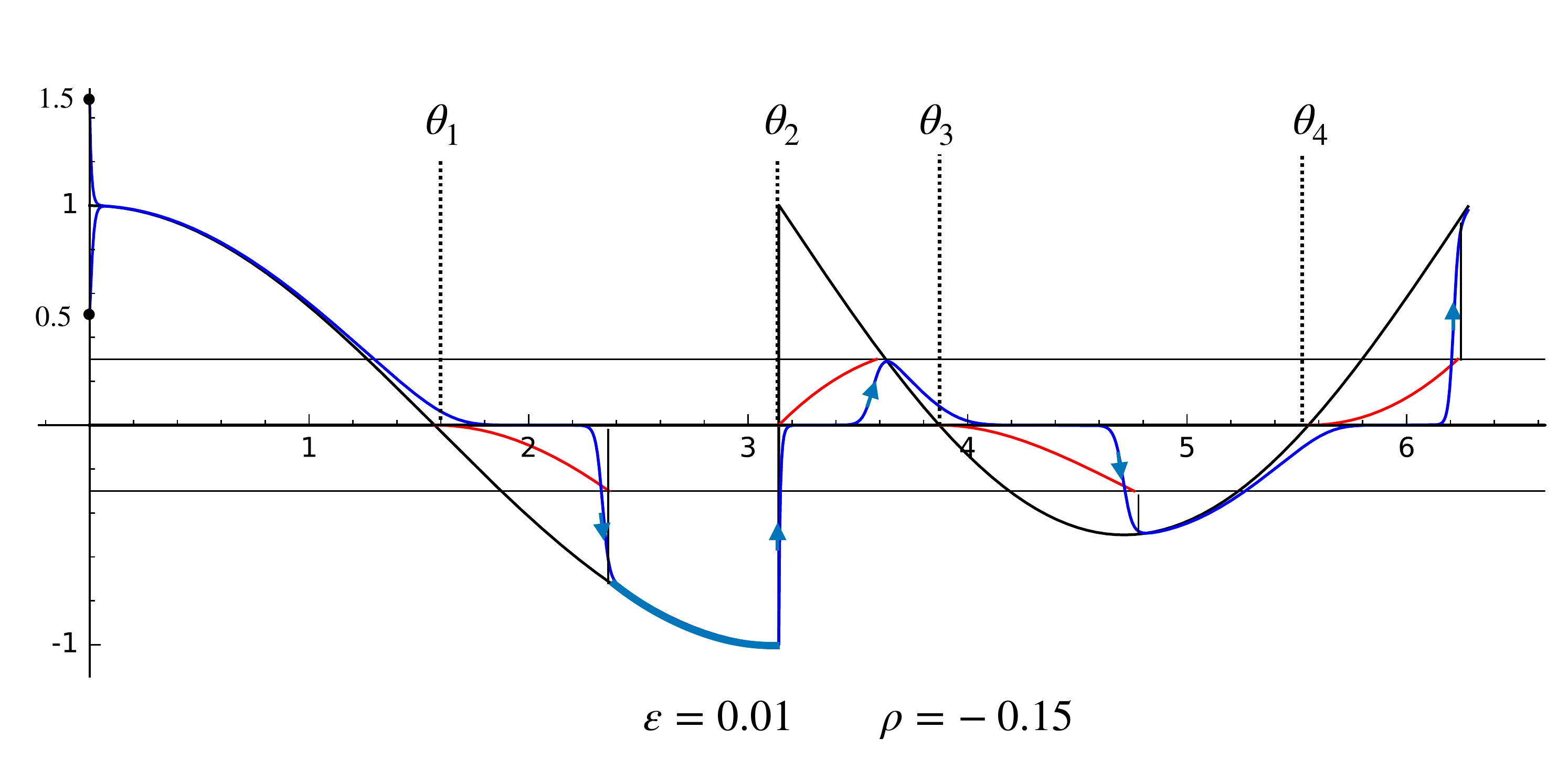}
 \caption{Two trajectories of $S_m$ that flank the periodic C-trajectory} \label{grapher1r2bister}
 \end{center}
 \end{figure}
The figure \ref{grapher1r2bis} represents the C-trajectory coming from the point: $(0,1)$. It is periodic. On the figure \ref{grapher1r2bister} we have the solutions from the initial conditions $(0,0.5)$ and $(0,1.5)$ respectively. We can apply the approximation theorem to them and thus, after having gone along respectively the segments $[(0,5), (0,1)]$ and $[(0,1.5), (0,1)]$ these two trajectories go along the periodic C-trajectory; this shows that the image of the segment $[(0. 5),(0,1.5)]$ by the Poincaré application is contained in the halo of $(0,1)$ so that the periodic solution of $S_m$ is infinitely close to the periodic C-path. Since $\max{S(\theta_i) - \theta_i}$ tends to $0$ when $|\rho|$ tends to $0$ this completes the proof of the proposition \ref{conjecturekat}.

 \section*{Acknowledgements}
 I warmly thank M. Benaïm, T. Sari and E. Strickler for the many good mood discussions during a (too !) short stay at the University of Neuchâtel at the end of the summer 2021.

\newpage
{\huge\textbf{Appendix}}
 \appendix
 \section{Introduction to NSA}\label{ANS}
 
The invention of NSA is credited to A. Robinson \cite{ROB16}. There are many versions of it, and the review paper \cite{FLE17} can be consulted to get an idea. The version on which we rely here, called {\em Internal Set Theory}, (I.S.T.) is due to E. Nelson. He proposed it in the article \cite{NEL77} which also contains an introduction to the practice of I.S.T. More than forty years after the publication of this article, it must be said that the hopes of Nelson and the other supporters of NSA to see this practice become widespread have been disappointed. Perhaps the price of the "entry ticket" was too high.

Afterwards Nelson also proposed to use more elementary versions than I. S. T., certainly less powerful, but whose fruitfulness he nevertheless showed in his book {\em Radically Elementary Probability Theory} \cite{NEL87}. Following the example of what his friend Reeb proposed (\cite{REE79}) (see also \cite{LUT87,LUT92} ) he constantly encouraged mathematicians to use weak theories of NSA as for example in this extract from his article {\em The virtue of simplicity} \cite{NEL07}. 
\dcom
Much of mathematics is intrinsically complex, and there is delight to be found in mastering complexity. But there can also be an extrinsic complexity arising from unnecessarily complicated ways of expressing intuitive mathematical ideas. Heretofore nonstandard analysis has been used primarily to simplify proofs of theorems. But it can also be used to simplify theories. There are several reasons for doing this. First and foremost is the aesthetic impulse, to create beauty. Second and very important is our obligation to the larger scientific community, to make our theories more accessible to those who need to use them. To simplify theories we need to {\em have the courage to leave results in simple, external form}\footnote{my emphasis}-fully to embrace nonstandard analysis as a new paradigm for mathematics.
\fcom
I have tried to comply with this recommendation, which I discuss a little later in the appendix.

 \subsection{To be standard or not, internal or not: the extended language.}
 For those who do not have elementary notions of logic, which is the case of many mathematicians, the practice of Nelson's I.S.T. theory is a bit confusing at first. Let us recall some elementary notions of logic without which one cannot understand the approach of I.S.T.
 
  The {\em formal set theory} Z.F.C. (see the classic \cite{KRI69} for a very accessible presentation) is the theory to which almost all contemporary mathematical texts are supposed to refer. 
  Let us recall that in this theory the formulas include in addition to the logical symbols ($\forall,\exists,\ =,\cdots$ ) a binary predicate $x \in y$ which reads ''$x$ belongs to $y$''. Z.F.C. proposes a series of axioms which allow us to construct from a primitive set, the empty set, new sets which have familiar names : 
  $$\Nmat, \mathbb{Z}, \mathbb{Q}, \Rmat, \cdots $$ 
  
  Functions are also considered to be sets. For example:
 $$(x \mapsto x^2 )\stackrel{\mathrm{def}}{\longleftrightarrow} \{(x,y)\in \Rmat^2\;;\; y = x^2\}$$
 In the expression above, after the ";", inside the  brackets $\{\}$, we have the writing "$y = x^2$". This is a  "formula" (a "statement") $\phi(x,y)$ which is or is not satisfied by the couple $(x,y)$. The fact that $ \{(x,y)\in \Rmat^2 \;;\; y = x^2\}$ can be considered as a set is an axiom, precisely the comprehension axiom :
 \bitbul
 \item $A_c \quad \quad$ S i $\phi$ is a ''formula'':$ \quad \forall x; \exists y; \forall z \;:\; (z \in y \Leftrightarrow z \in y \;\mathrm{et} \phi(z))$
 \fit

 A "formula" is not just anything. It is a sequence of parentheses, logical symbols, occurrences of the membership symbol $\in$ and constants (the names of the sets already constructed) which respects a construction syntax ; correctly written formulas constitute the $\mathcal{L}$ language of Z.F.C. A given formula $\phi$  can be true or false.
 
\paragraph{ I.S.T. is an extension of Z.F.C.} We increase the possibilities of expression of $\mathcal{L}$ by adding the predicate $st(x)$ which reads "$x$ {\em is standard}" ; the rules of manipulation of $st$ are defined by three axioms, {\em Idealization}, {\em Standardization} and {\em Transfert}  which I do not specify at once\footnote{Except specialists of set theory, which mathematicians have in mind, when they write their demonstrations, the ten or so axioms of Z.F.C.?} from which useful properties of $st$ can be quickly deduced. For example we have :
 \ben
 \item $st(0) \longleftrightarrow$ $0$ is standard.
 \item $n \in \Nmat : st(n) \Rightarrow st(n+1) \longleftrightarrow$ the successor of a standard is standard.
 \item $\exists \; \omega \in \Nmat \; \forall n \in \Nmat : (st(n) \Rightarrow n < \omega) \longleftrightarrow$ there exists an infinitely large integer (greater than all the standard integers).
 \fen
 One may be tempted to see a problem here: a classical theorem of Z.F.C. says that if $E$ is a part of $\Nmat$, if $0\in E$ and if $n \in E \Rightarrow (n+1) \in E$ then $E = \Nmat$. Don't the points 1. + 2. above imply that the set of standard integers is the whole  $\Nmat$ ? This contradicts 3. In reality there is no problem because nothing allows us to assert that the standard integers constitute a set: $st(x)$ not being a formula of the language $\mathcal{L}$ of Z.F.C. the axiom $A_c$ cannot be invoked.

 That's why, from now on, we must have in mind that there are two kinds of formulas;
 \bito 
 \item The  internal formulas $\phi$, those belonging to the $\mathcal{L}$ language of Z.F.C. for which we can write, thanks to the  comprehension axiom :
 $$ E = \{x \in y \;;\; \phi(x)\}$$ 
 where $E$ is a set of the formal theory, 
 \item the  formulas $\Phi$ that are not internal for which the writing :
 $$ \{x \in y \;;\; \Phi(x)\}$$ 
 does not designate a set; to avoid the risk of believing that we are dealing with a set, when $\Phi$ is not an internal formula we can write : 
 \beq \label{ensembleexterne}
\eset{x \in y \;;\; \Phi(x)}
 \feq
 \fit
 For example, in $\Rmat^2$ if $\mathcal{C}$ is the graph of a function $f$ we could write :
 $$ \mathrm{Halo}(\mathcal{C}) = \eset{ (x,y) \,:\, y \sim f(x)}$$ 
 (Here the formula $\Phi$, which deals with the pairs $(x,y)$, is the external formula $y \sim f(x)\;$).
 
 We must be careful to distinguish :
 \bitbul
 \item To be standard which is a quality that a set (a mathematical object) may or may not have.
 \item To be internal  which is a quality that a formula can have or not.
 \fit
 
 It can be demonstrated  that I.S.T. is a {\em conservative extension} of Z.F.C. which means that if A and B are two internal propositions (which are expressed in $\mathcal{L}$) and if there exists a proof in I.S.T. of $A \Rightarrow B$ then there exists a proof in $\mathcal{L}$. 
This last property has been misinterpreted as meaning that nothing new can be proved with I.S.T. This is true only for internal theorems, i.e. the implications $A \Rightarrow B$ where $A$ and $B$ are internal, but theorems which are not internal can be proved.
 The importance of conservativity lies in the fact that if I.S.T. was contradictory then Z.F.C. would be too. Indeed, if I.S.T. were contradictory there would be a proof (for example) of the internal proposition $1 = 2$, and, since I.S.T. is conservative, there would be a proof in Z.F.C. of $1 = 2$, hence a contradiction in Z.F.C. From the point of view of the risk of contradiction, hence of rigor, it is no more dangerous to work in I.S.T. than in Z.F.C.
 
  \subsection{Who is standard?}
  The axioms of I.S.T. are the rules for handling the predicate $st(x)$. They are as follows:
  We use the abbreviations\\
$\forall^{st}x \longleftrightarrow \forall x \;(x\; \mathrm{standard}) \Rightarrow$\\
$\forall^{fin}x \longleftrightarrow \forall x \;(x\; \mathrm{fin}) \Rightarrow \quad$\\
$\forall ^{st\;fin}\longleftrightarrow \forall^{st} x \;(x; \mathrm{fin}) \Rightarrow $\\
$\exists^{st} x \longleftrightarrow \exists x\; (x\; \mathrm{standard})\; \wedge\;$\\
The word "finite" is taken in its usual mathematical sense: a set $x$ is "finite" if any injection of $x$ onto itself is surjective.
\bitbul
\item {\em \textbf{I}dealization}~: Let $B(x,y)$ be an {\em internal} formula:
$$\forall ^{st\;fin} z\, \exists\, x \forall \,y \in z \;\; B(x,y) \Longleftrightarrow \exists\,x \;\forall^{st}y\; B(x,y)$$
\item {\em \textbf{S}tandardisation}~: Let $\phi(z)$ be a formula, internal or not:: 
$$\forall^{st}x\,\exists^{st}y\,\forall^{st}z\;(\;z \in y \Leftrightarrow z\in x \wedge \phi(z))$$
\item {\em \textbf{T}ransfer }: Let $A(x,t_1,\cdots,t_k)$ be an internal formula with no free variables other than $(x,t_1,\cdots,t_k)$. Then~:
$$\forall^{st}t_1 \cdots \forall^{st}t_k \; \big( \forall^{st}x\; A(x,t_1,\cdots,t_k) \Rightarrow \forall x \;A(x,t_1,\cdots,t_k) \big)$$
and by contraposition~:
$$\forall ^{st}t_1 \cdots \forall ^{st}t_k \; \big( \exists \; x;A(x,t_1,\cdots,t_k) \Rightarrow \exists ^{st} x \;A(x,t_1,\cdots,t_k) \big)$$
\fit
hence the name given to the theory: \textbf{I.S.T.}, which also stands for {\em Internal Set Theory}.

Note the similarity between the $A_c$ axiom of Z.F.C. and the axiom of standardisation. The former tells us that for any formula, the elements of a set $x$ which satisfy this formula constitute a set, the latter specifies that if the starting set is standard, the set whose standard elements satisfy the formula can be chosen as standard. As essentially the new sets of Z.F.C. are constructed by the application of $A_c$ it follows that all mathematical objects (the ''sets'') constructed in a unique way are standard.  Thus all the objects of traditional mathematics such as $\Nmat, \Rmat, x\mapsto \cos (x), \cdots$ are standard. But all the elements of a standard set are not necessarily standard. In fact all the elements of a standard set are standard if and only if the latter is a finite set. All this leads to the following situation which at first sight seems strange: 
  \ben
  \item $\Nmat$ is standard.
  \item $0$ is standard.
  \item $st(n) \Rightarrow st(n+1)$
  \item $\exists \omega \in \Nmat$ and $\omega$ is not standard.
  \fen
On the one hand the successive ordinals $0,1,2, \cdots$ are constructed one after the other, they are standard. On the other hand the existence of $\Nmat$ is deduced almost directly from the axioms of Z.F.C. (There exists an infinite set, there exists a good order on any set, etc.) so it is standard. But in the universe of Z.F.C. there is nothing to say that the collection of ordinals is equal to the set $\Nmat$
 which leaves room for an element of $\Nmat$ which would not be an ordinal constructed step by step,i.e. a nonstandard integer: we therefore add the axiom that among all the elements of $\Nmat$ there is one which is not standard.

A few examples shed light on the situation.
\bit
\item $\Rmat$ is standard, it contains non-standard elements, the infinitesimals, the infinitely large, but not only. For example $1+\eps$, with $\eps$ infinitesimal, is a non-standard, neither infinite nor infinitesimal.
\item $\Rmat^2$, $\Rmat^3$ are standard, $\Rmat^n$ is standard if $n$ is standard.
\fit
 The function of $\Rmat^2$ in $\Rmat$ :
$$(x, m) \mapsto \frac{\emat^{mx}-\emat^{-mx} }{\emat^{mx}+\emat^{-mx}} = \tanh(mx)$$ 
 is a "true" function in the traditional sense; whatever $m$ and $x$ are, whether they are standard or not, the real $\emat^{mx}$ is well defined internally by formulae (here the sum of a convergent series) which do not use the predicate $st()$. It is a standard function of the couple $m, x$. But if now $m$ is understood as a parameter, for all $m$ in $\Rmat$ :
 $$x \mapsto f_m(x) = \mathrm{th}(mx)$$
 is a function of $\Rmat$ in $\Rmat$ which, as a function, is standard if $m$ is standard and non-standard if $m$ is not. Take $m$ i.g., so $\mathrm{th}(mx)$ is non-standard, but it has not ceased to be the traditional function $\mathrm{th}(mx)$ which is continuous at the point $x = 0$. However as soon as $x$ is not infinitesimal, the product $mx$ is infinitely large (negative or positive) and the value taken by $mathrm{th}(mx)$ is infinitelsimaly close to $±1$. When $m$ is infinitely large the function $f_m$ has the particularity that an infinitesimal increase of the variable $x$ in the neighbourhood of $0$ can lead to a non-infinitesimal increase of $f$. We say that the function is ''continuous'' (in the traditional sense) but is not S-continuous (i.e. a function $f$ is said to be S-continuous if $dx \sim 0$ limplies to $f(x+dx) \sim f(x)$).
 
  All this is just a way of talking in the new extended language about the (non-uniform) convergence of the family of functions $f_m$ to the function $\mathrm{sign}(x)$ when $m$ tends to infinity. For this example the benefit
 does not seem obvious. But it becomes so in the case of singular perturbations of differential equations.
 
 Indeed, in this case, we are interested in a family of  "objects" called ''differential equations":
 \beq
\Sigma_{\eps} \quad\quad \quad \displaystyle \eps \frac{dx}{dt} = f(x) \Longleftrightarrow \frac{dx}{dt} = \frac{1}{\eps} f(x)
 \feq 
  when $\eps$ tends to $0$. One would like to prove theorems of the kind ''the phase portait of $\Sigma_m$ tends towards that of $\Sigma_0$'' but one cannot because $$0\, \frac{dx}{dt} = f(x)$$ is not a differential equation. In NSA we fix $\eps >0$ infinitely small and we are interested in the phase portrait of $\Sigma_{\eps}$ which remains perfectly defined as the set of trajectories of $\Sigma_{\eps}$ and we will seek to establish results independent of the value of $\eps$, provided that the latter is infinitely small as it is in the trajectory approximation theorem of this article.
  
   A last word on the standard functions that interest us here. Let $x \mapsto f(x)$ be a standard real-valued function whose lower bound $\inf(f)$ is strictly positive; then $\inf(f) \gnsim 0$. Indeed, the lower bound being uniquely defined from a standard object is, in its turn, standard; a strictly positive standard number is not infinitesimal.

\subsection{Weaker theories than I.S.T.}\label{cheapANS}

It is undeniable that the formulas that define the I.S. and T. axioms have a repulsive effect on the ordinary mathematician. This is why mathematicians like G. Reeb and R. Lutz have proposed to simplify this axiomatic by weakening it. In \cite{LUT87} a ''sub-system'' of I.S.T. called Z.F.L. (''L.'' in honour of Leibnitz) is proposed, of which Lutz tells us in \cite{LUT92} :
\dcom
In I.S.T. there are automatic equivalences between internal theorems and external formulations (see Nelson \cite{NEL77}). This is not the case in Z.F.L.; competing formulations of the intuitive ideas of limit, continuity etc. coexist while remaining irreducible to each other. 
\fcom
Nelson himself proposed a weakened version of I.S.T. in his book {\em Radically Elementary Probability Theory} where he introduces an appendix showing the formal equivalence of the radically elementary results with the traditional results by the remarks :
\dcom
The purpose of this appendix is to demonstrate that theorems of the conventional theory of stochastic processes can be derived from their elementary analogues by arguments of the type usually described as generalized nonsenses; there is no probabilistic reasoning in this appendix. This shows that the elementary nonstandard theory of stochastic processes can be used to derive conventional results; on the other and, it shows that neither the elaborate machinery of the conventional theory nor the devices from the full theory of nonstandard analysis, needed to prove the equivalence of the elementary results with their conventional forms, add anything of significance : the elementary theory has the same scientific content as the conventional theory. This is intended as a self-destructing appendix.
\fcom

These authors propose a simplified external language in which only the structure of real numbers remains, enriched by the possibility of talking about orders of magnitude, and they encourage us to work on it; the price of this simplicity is the impossibility of translating into conventional language what is said externally. But is this really a price to pay? Is it not, on the contrary, more rewarding to keep both points of view? We have noticed that for $m$ infinitely large the function $m \mapsto \tanh (mx)$ is (traditionally) continuous, but not {\em S-continuous}. On the other hand the function $x \mapsto \eps [x/ \eps]$ (where $[\cdot]$ denotes the integer part), when $\eps \sim 0$, is S-continuous, but is not continuous: it progresses by small infinitesimal jumps and does not satisfy the intermediate value theorem. Continuity and S-continuity do not speak exactly the same thing when considered on non-standard objects. 

On the other hand this question of the return of external statements to internal statements of $Z.F.C.$ is related to the use of the axiom of choice. A criticism often formulated against the use of the NSA is its essentially ''non-constructive'' character. It is true that the demonstration of the existence of a NSA model with a strength at least equal to I.S.T. inside $Z.F.C.$ is based on the axiom of choice in its most general version. But this is not true for weak versions. Logicians have shown that, provided that we give up transforming certain external statements into internal ones, suitable weak axiomatic versions allow 
a nonstantard practice of the essentials of mathematics by using, at most, only the axiom of countable choice (see \cite{HRB21} for a relatively accessible account).

It is for all these reasons that I have tried not to use all the power of I.S.T., in particular the axiom $S.$ which allows to define {\em the shadow} of a function as the {\em standardized} of the {\em halo} of its graph, which would have obliged the reader to understand what the standardized is, and therefore would have required an additional effort. 
On the other hand, I have not tried to define precisely the non-standard axiomatic system that would best formalize the natural intuition of practitioners of the theory of singularly perturbed differential equations. This is a work that remains to be done.

 \section{ANS and dynamical systems} \label{ANSandEDO}

 \subsection{Continuous dependence of solutions} \label{gronwall}

 The very classical theorem on continuous dependence of the solutions of an ODE on the initial conditions and on the perturbations that we deduce from Gronwall's lemma has the following external expression :
 \begin{theoreme}
 Let $ f $ be $C^1$-limited and globally lipschitz,with limited (not infinitely large)  constant $k$. Let $g(x)$ be of class $C^1$
infinitesimal (i.e. for all $x$ we have $g(x) \sim 0$). Let $x_0\sim y_0$ and let $x(t)$ and $y(t)$ be the respective solutions of :\\
 $\frac{dx}{dt} = f(x)\quad x(0) = x_0$\\
$\frac{dy}{dt} = f(y)+g(y);\quad y(0) = y_0$\\
 Then, for all limited $T$ :
 $$t \in [0,T] \Longrightarrow y(t) \sim x(t)$$
 \end{theoreme}

 \textbf{Proof} Let $E$ be a compact set containing the trajectory $\{ x(t); t \in [0, T]\}$ in its interior. The maximum of $|g(x)|$ on $E$ is $m \sim 0$. The classical Gronwall inequality says that, as long as $y(t) \in E$ we have :
 $$ |x(t)-y(t)| \leq |x_0-y_0|\emat^{kt m} \sim 0$$
as $k$, $t$ and $m$ are limited so is $\emat^{kt m}$. So $y(t)$ cannot exit from $E$ before $T$ and $y(t) \sim x(t)$ until $T$.
 $\Box$\\
 We leave it to the reader to formulate the vector case and to get rid of the assumption of the existence of a global Lipschitz constant for $f$.

 \subsection{Entering the {halo} of a slow curve.}\label{halocourbelente}
 Let's start by recalling a well-known elementary result concerning vector fields in the plane. Consider a  curve $\Gamma$, closed, with no intersection, \\
 \begin{minipage}{0.55 \textwidth} 
 piecewise differentiable, consisting of the union of segments $\arc{ab}$, $\arc{bc}\cdots \arc{ga}$ (figure on the right), which defines a domain $\Gamma$ and a differentiable vector field $X$, {\em which does not vanish} in $\Gamma \cup \partial \Gamma$ and which points strictly to the interior of $\Gamma$, except for one of the segments ($\arc{ef}$ on the figure) where it points to the exterior. 
 \end{minipage}
 \begin{minipage}{0.40 \textwidth}
 \includegraphics[width=1\textwidth]{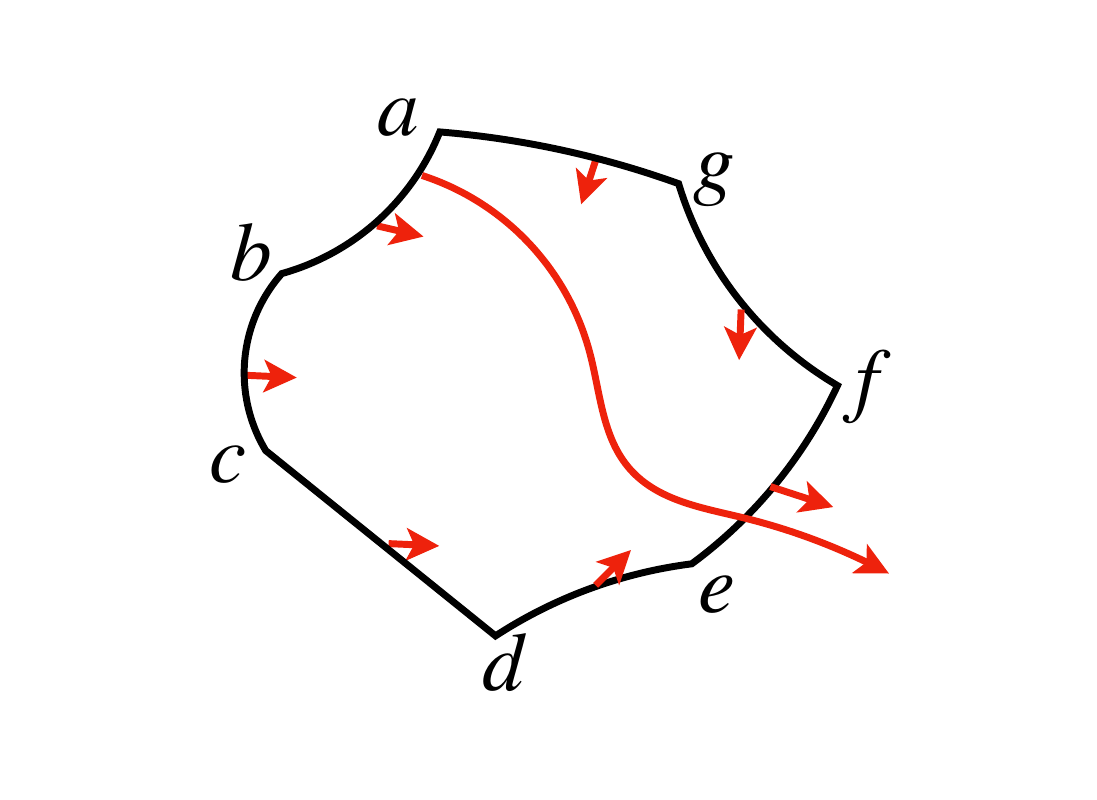}
 \end{minipage}
 \begin{theoreme}\label{piege}
 For any initial condition in $\Gamma \cup \partial \Gamma$ the trajectory the corresponding of $X$  exits of $\Gamma$ at some point of $\arc{ef}$.
 \end{theoreme}
  
  We now consider the system :
 \beq \label{eqLR3}
\quad \quad \quad  \left\{
\begin{array}{lcl}
\displaystyle \frac{dx}{dt}& =& 1\\[6pt]
\displaystyle \frac{dy}{dt} &=& \displaystyle \frac{1}{\eps}\Big(f(x)-y\Big)
 \end{array} 
 \right.
\feq
with $f$ $C^1$-limited. We suppose $\eps \sim 0$ fixed. 0n says that such a system is  a{\em slow-fast} system. The graph of $f$ is called the {\em slow curve}.

Let $(x_0,y_0)$ be a limited initial condition at time $0$, not belonging to the halo $\mathcal{H}$ of the slow curve (above the graph of $f$ to fix ideas).
\begin{proposition}\label{cl1}
Let $(x(t),y(t))$ be the solution of \eqref{eqLR3} with initial condition $(x_0,y_0)$. There exists $t^*>0$ :
\ben
\item $t^*\sim 0$
\item $t \in \,[0,t^*] \Longrightarrow x(t) \sim x_0$
\item $ y(t^*) \sim f(x_0)$
\item $ t > t*$ and limited $ \Longrightarrow y(t) \sim f(x(t))$
\fen
In other more geometrical terms: the trajectory goes along the segment $\{(x_0,y) : y \in [f(x_0),y_0]\}$, enters the halo of the graph of $f$ before the time $t^* \sim 0$ and stays there as long as $t$ is bounded (beyond that we can say nothing as shown by the example $f(x) = x^2$ which integrates explicitly). 
\end{proposition}
\textbf{Proof.}
We suppose, to fix the ideas, that $f'(x_0) < 0$.
 Let $\alpha \gnsim 0$ be such that the graph of $f(x)+\alpha$ is below $y_0$; there are some, since $(x_0,y_0)$ is not in the halo of the graph of $f$. Let $\theta \gnsim 0$. Consider the closed path $a,b,c,d$ defined by :\\
  \begin{minipage}{0.50 \textwidth} 
  \bito 
 \item Segment $\arc{a,b}$ : line segment joining $a$ to $b$ ; $a$ is the initial condition $(x_0,y_0)$, $b$ is the point $(x_0,f(x_0))$ of intersection of the vertical $x = x_0$ with the graph of $f$.
  \item Segment $\arc{b,c}$: portion of the graph of $f$ between $b$ and $c$; $c$ is the intersection point of the graph of $f$ with the vertical $x = x_0+\theta$.
  \fit
\end{minipage}
  \begin{minipage}{0.5 \textwidth}
 \includegraphics[width=1\textwidth]{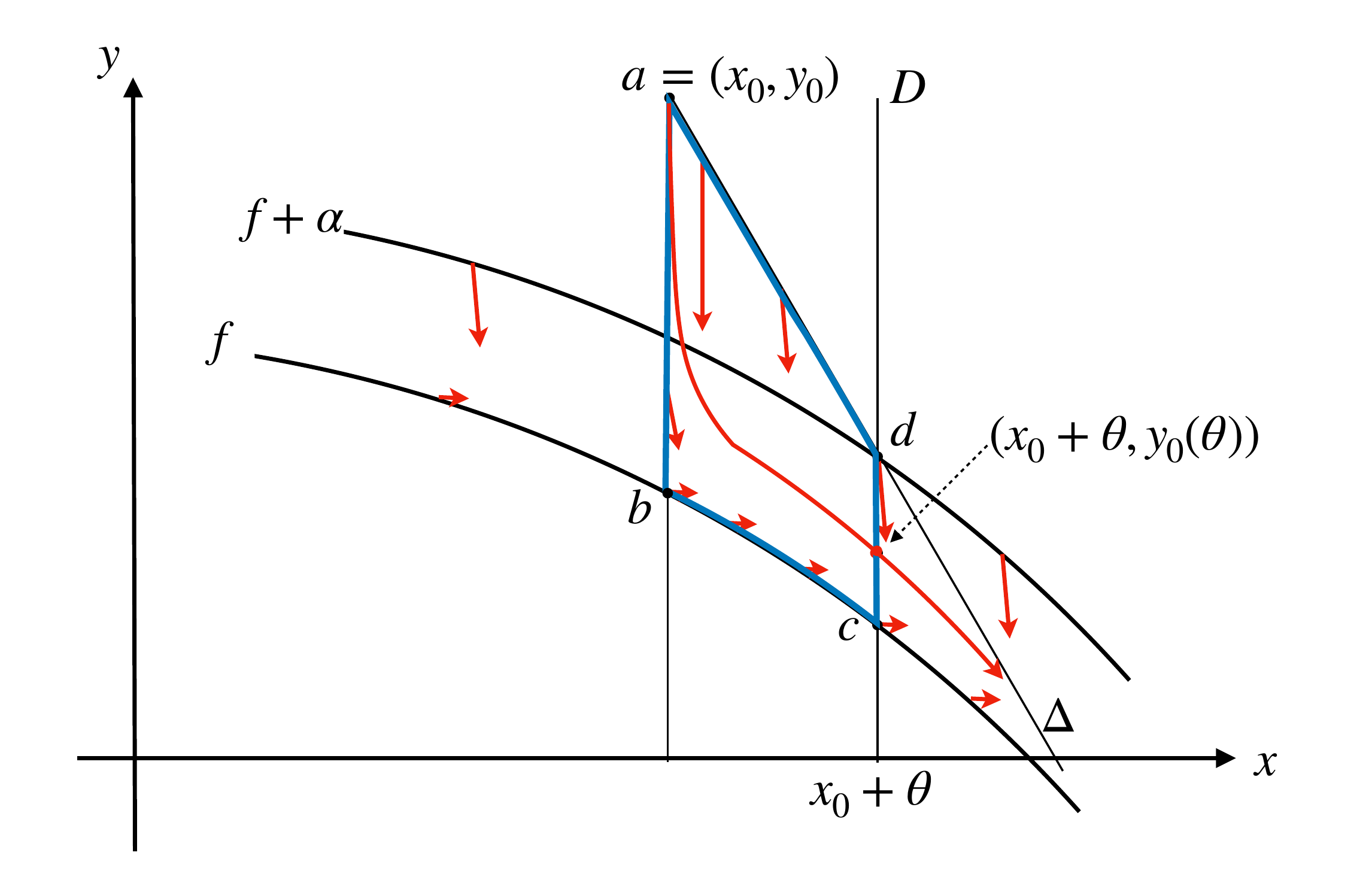}
 \end{minipage}
 \bito 
  \item Segment $\arc{c,d}$ : line segment joining $c$ to $d$ ; $d$ is the intersection point of the graph of $f+\alpha$ with the vertical $x = x_0+\theta$.
 \item Segment $\arc{da}$: line segment joining $d$ to $a$.
 \fit 
On all these segments, except the segment $\arc{c,d}$, the vector field defined by \eqref{eqLR3}) points strictly to the interior of $\Gamma$ ; for the segment $\arc{d,a}$ this is a consequence of the fact that for $(x,y)$ above the graph of $f(x)+\alpha$ we have $\frac{dy}{dx} < \frac{ f(x)-y}{\eps}$ which is {infinitely large} (negative), so smaller than the slope of $\Delta$. which is not infinitely large because $\theta \nsim 0$.
This shows that for all $\alpha \gnsim 0$ we have $f(x_0) < y(\theta) < f(x_0)$ so that $y(\theta)$ is infinitely close to $f(x_0)$, so that $(x(\theta),y(\theta))$ is in the {halo} of the graph of $f$. We leave it to the reader to make the obvious modifications that allow us to obtain the same result in the case where the derivative of $f$ in $x_0$ is positive or zero.

We have shown that for any $t^*\sim 0$ the solution $(x(t), y(t)) $ is in the {halo} of the graph of $f$, but the proposition says a little more: that there exists a $t^*\sim 0$ such that $\forall t \,\geq t^* \;:\; y(t)\sim f(y(t))$. This last point is obtained by means of a typically ''non-standard'' argument (called Robinson's lemma) which is the following.
Consider the set :
$$
E =\left \{ t \in [0,1] \;:\; \frac{ |y(t)-f(x(t))|}{t} \leq 1 \right \}
$$
The definition of $E$ does not involve any term of the new language, so it is a ''true'' set in the usual set theory. This set contains all the real numbers not {infinitesimals} of $[0;1]$ ; indeed, for $t\nsim 0$ we have $|y(t)-f(x(t))| \sim 0$ which means that $\frac{ |y(t)-f(x(t))|}{t} \leq 1$ (we divide an infinitesimal by a non {infinitesimal}, the result is {infinitesimal}).    The lower bound $l$ of $E$ is infinitesimal because if it were not, $l/2$ would not be infinitesimal and so $l/2$                          
would be in $E$ which contradicts the definition of the lower bound. If the lower bound of $E$ is infinitesimal, $E$ necessarily contains an infinitesimal $t^*$. Since $t^*$ is in $E$ we have :
$$ \frac{|y(t^*)-f(x(t^*))|}{t^*} < 1$$
For the quotient to be smaller than 1 the numerator must be {infinitesimal} since the denominator is.\\
$\Box$\\\\
In the previous result it is essential that $|f'(x)|$ is limited. If this were not the case, the point $c$ could be rejected in the negative infinity and we could no longer guarantee that the segment $\arc{da}$ has a limited slope, so that, along $\arc{da}$, the field points to the interior of the domain. 

With obvious adaptations of this proof one obtains the following result a little more general than \eqref{eqLR3} which can be considered as the NSA version of Tychonov's theorem (in dimension 2) (see \cite{TYK52, LOB98}) :\\
Let the system be :
 \beq \label{eqLR4}
 \left\{
\begin{array}{lcl}
\displaystyle \frac{dx}{dt}& =& h(x,y)\\[6pt]
\displaystyle \frac{dy}{dt} &=& \displaystyle \frac{1}{\eps}g(x,y)\Big(f(x)-y\Big)
 \end{array} 
 \right.
\feq
where $h$ is a $C^1$-limited function that does not vanish.
\begin{proposition}\label{cl2}
Consider the system \eqref{eqLR4}. We assume that:
$$\min g(x,y) \gnsim 0$$ 
Let $(x(t),y(t))$ be the solution of \eqref{eqLR3} with initial condition $(x_0,y_0)$. There exists $t^*>0 $ :
\ben
\item $t^* \sim 0$
\item $t \in \,[0,t^*] \Longrightarrow x(t) \sim x_0$
\item $ y(t^*) \sim f(x_0)$
\item $ t > t*$ {\em and limited} $ \Longrightarrow y(t) \sim f(x(t))$
\fen
\end{proposition}

\subsection{Semi-slow-fast systems.}\label{demiLR}
A semi-slow-fast field is a field which is infinitely large on one side of a slow curve, not infinitely large on the other side. I deal with the particular case of the semi-slow-fast field that concerns us. 

Let be the vector fields :
\beq \label{DLR1}
\Sigma \quad \left \{
\begin{array}{lcl}
\displaystyle \frac{dx}{dt}& =&1 \\[6pt]
\displaystyle \frac{dz}{dt} & =&\displaystyle z\sqrt{1+\varphi (z) } ( f(x)-\psi(z))
\end{array} 
\right.
\feq
and
\beq \label{DLR2}
\tilde{\Sigma} \quad \left \{
\begin{array}{lcl}
\displaystyle \frac{dx}{dt}& =&1 \\[6pt]
\displaystyle \frac{dz}{dt} & =&\displaystyle zf(x)
\end{array} 
\right. \quad \quad \quad \quad \quad 
\feq
We immediately verify that under the following hypotheses :
\ben
\item $f$ does not cancel is $C^1$-limited and $inf(f(x)) \gnsim 0$
\item $\varphi > 0$ and $psi >0$
\item $z \gnsim 1 \Longrightarrow \psi(z) \;\mathrm{and}$ {i.g.}
\item $z \lnsim 1 \Longrightarrow \varphi(z)\sim 0 \;\mathrm{et}\; \psi(z)\sim 0$
\fen
one has :
\bitbul 
\item$z \lnsim 1 \Longrightarrow z\sqrt{1+varphi (z) } ( f(x)-\psi(z)) \sim zf(x) \quad \quad z \gnsim 1 \Longrightarrow \frac{dz}{dt} = - \mathrm{i.g.}$  
\fit 
So, below the line $z = 1$ we have $\Sigma \sim \tilde{\Sigma}$\\\\
\textbf{The attractive case: $f>0$}. The same argument as in the slow-fast case (appendix \ref{halocourbelente}) allows us to state that if $(x(t),y(t))$ is the trajectory coming from the point $a$ (above the line $z =1$) there exists $t \sim0$ such that $z(t) \sim 1$.\\
\begin{minipage}{0.5 \textwidth}
Let us now start from the point $c$ located below $z = 1$. Let $(\tilde{x}(t),\tilde{y}(t)) $ be the solution of $\tilde{\Sigma}$ at time $0$ and $\tilde{y}(\tau_1) = 0$.  Since $inf(f(x)) \gnsim 0$  for any $\alpha \gnsim 0$ we have $\tilde{y}(\tau_1-\alpha) \lnsim 1$ and thus $\Sigma \sim \tilde{\Sigma}$; according to Gronwall's lemma (appendix \ref{gronwall}) we have therefore:
$$ 0 \lnsim \alpha \leq \tau_1 \Longrightarrow |y(t-\alpha)-\tilde{y}(t-\alpha)| \sim 0$$
The set :
\end{minipage}
\begin{minipage}{0.5 \textwidth}
\includegraphics[width=1\textwidth]{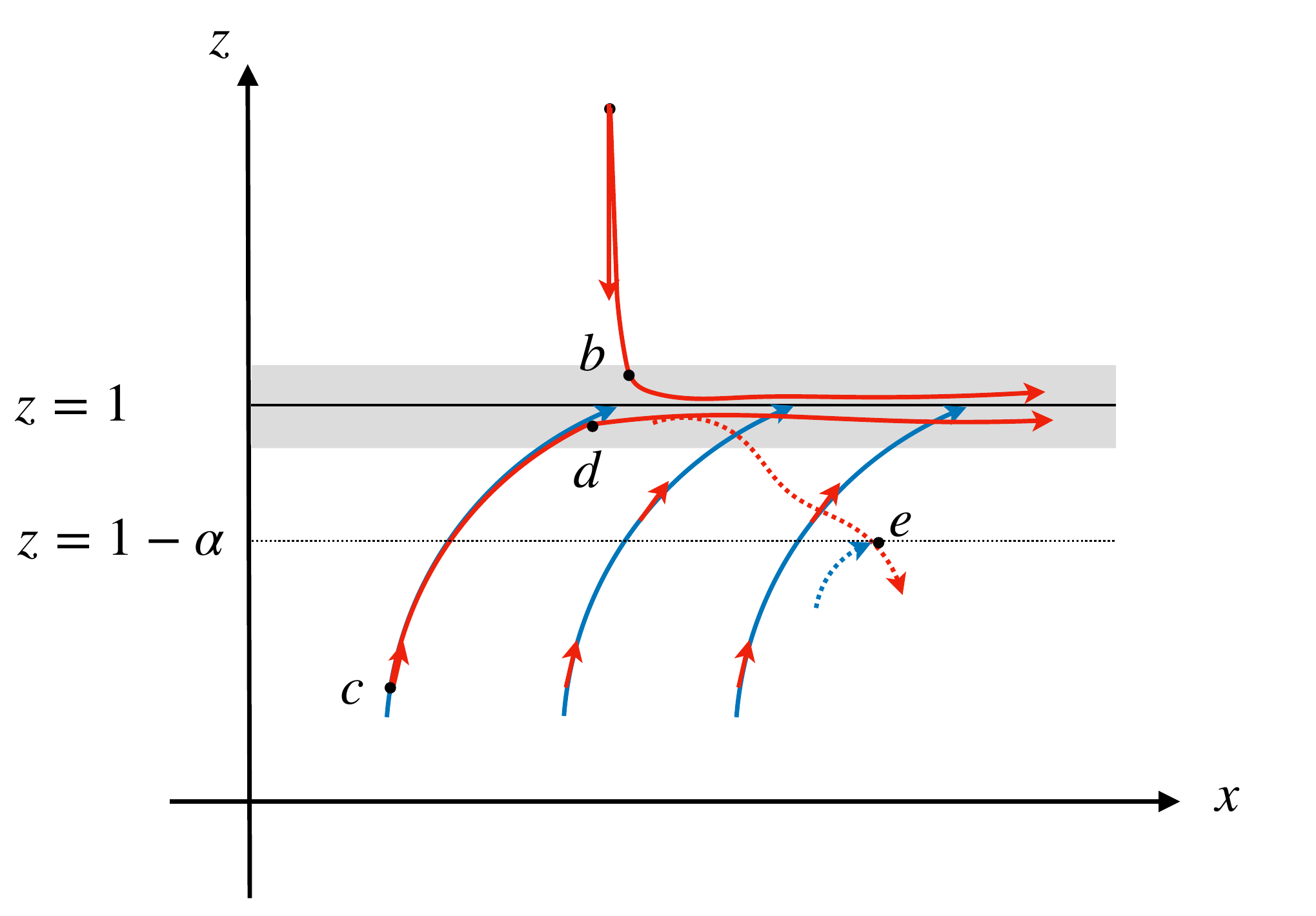}
\end{minipage}\\[6pt]
$$E = \left \{\alpha \in [0,\tau_1] : \frac{ |y(t-\alpha)-\tilde{y}(t-\alpha)|}{\alpha} < 1 \right \} $$
contains all the non infinitesimals, it therefore contains a $\alpha \sim 0$ which shows that there exists $t_1\sim \tau _1$ such that $(x(t_1),y(t_1)) \sim (\tilde{x}(t_1),\tilde{y}(t_1)) $. The trajectory of $\Sigma$ originating from $c$ at time $0$ enters the halo of $z = 1$ at a time $t_1$ infinitely close to the time $\tau_1$ where $(\tilde{x}(t),\tilde{y}(t)) $ reaches $z = 1$ and remains infinitely close to $(\tilde{x}(t),\tilde{y}(t)) $ on $[0, t_1]$. From this instant the trajectory remains infinitely close to $z = 1$ ; indeed for any $\alpha \nsim 0$ if $t_2$ is the first instant where $y(t) = 1-\alpha$ at the point $e = (x(t_2),y(t_2)$ we have a contradiction with the fact that at this point the field points strictly upwards.\\\\
\textbf{The crossing case: $f< 0$}.\\
\begin{minipage}{0.5 \textwidth}
Let the point $a =(x_0,y_0)$ above $z = 1$, the point $ c= (x_0,1)$ and $(\tilde{x}(t),\tilde{y}(t))$ be the trajectory of $\tilde{\Sigma}$. The trajectory of $\tilde{\Sigma}$ coming from $a$ at time $0$, enters the halo of $z = 1$ at time $t_1 \sim 0$, comes out under it at time $t_2 \sim 0$ and then remains infinitely close to the trajectory $(\tilde{x}(t),\tilde{y}(t))$. The duration of the crossing of the halo of $z = 1$ is infinitely small because $\frac{dz}{dt} < z f(x) \leq \min(zf(x)) \lnsim 0$.
\end{minipage} $\quad$
\begin{minipage}{0.45 \textwidth}
\includegraphics[width=1\textwidth]{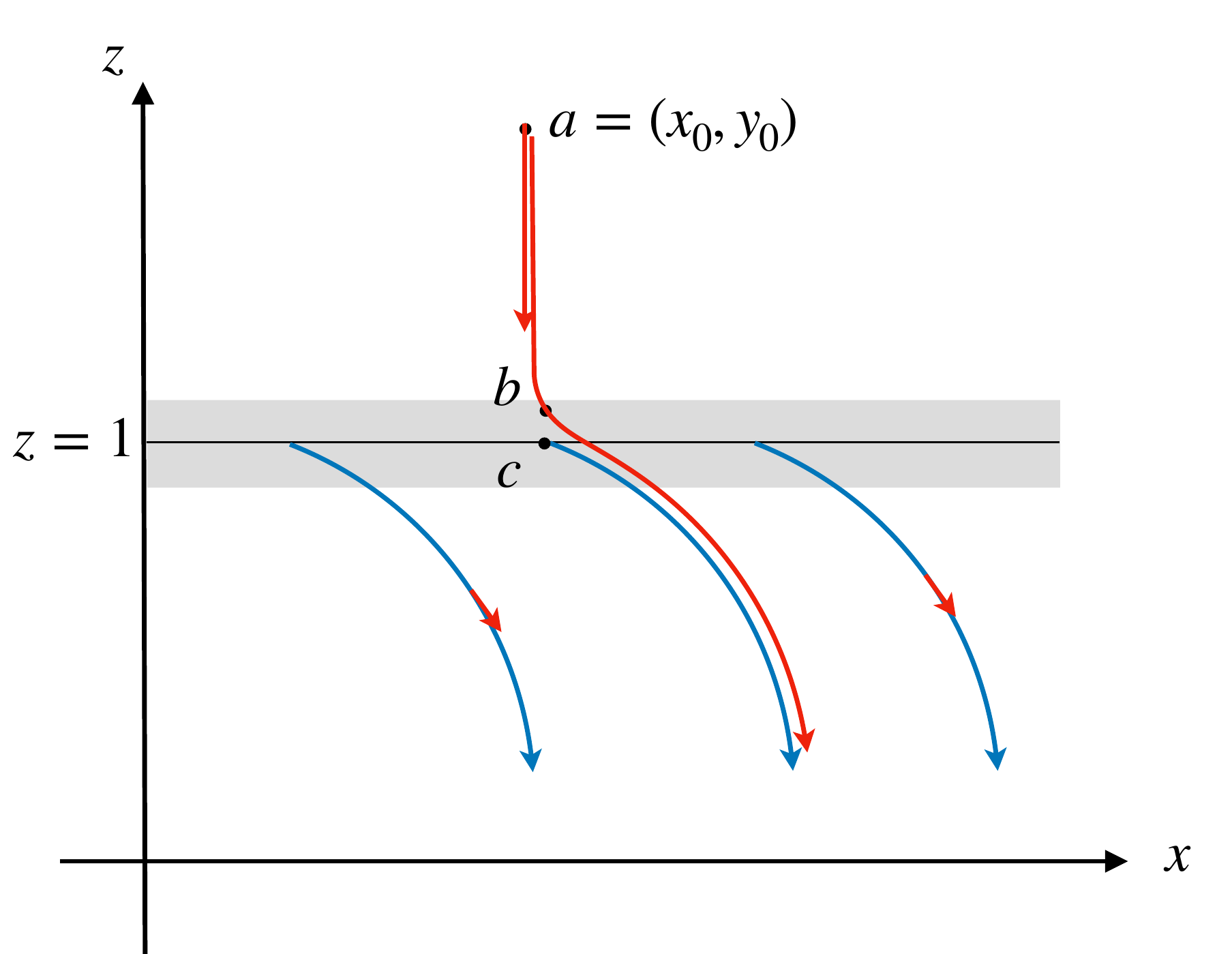}
\end{minipage}{6pt}

\subsection{Passing a change of sign of $f$.}{\label{changementdesigne}
 Let us consider for $z > 0,\; \rho \lnsim 0, \eps \sim 0$ the system :
 \beq \label{eqchsigne}
F \quad \quad \quad \left\{
\begin{array}{lcl}
\displaystyle \frac{dx}{dt}& =&1 \\[6pt]
\displaystyle \frac{dz}{dt} & =&\displaystyle z\sqrt{ 1+\exp \left(2 \frac{\rho-\ln(z)}{\eps} \right) } \left( f(x)-z^{\frac{1}{\eps}}\right)
\end{array} 
\right.
\feq
 and we make the following assumptions:
 \begin{hypothese} Assumptions about $f$:
 \bitbul
\item For $ x > 0$ the function $f$ is equal to a function $f^+$, $C^1$-limited defined on $\Rmat$ strictly negative for $x>0$
\item For $ x < 0$ the function $f$ is equal to a function $f^-$, $C^1$-limited  defined on $\Rmat$ strictly positive for $x<0$
 \fit
 \end{hypothese}
We introduce the vector field defined by :
  \beq \label{eqchsignealphaplus}
\tilde{F}  \quad \quad \left \{
\begin{array}{lcl}
\displaystyle \frac{dx}{dt}& =&1\\[6pt]
\displaystyle \frac{dz}{dt} & =&\displaystyle z f(z)
\end{array} 
\right.
\feq
 and we note respectively $\big(x(t,(x_0,z_0)),z(t,(x_0,z_0))$ and $\big(\tilde{x}(t,(x_0,z_0)),\tilde{z} (t, (x_0,z_0))$ the solutions of initial condition $(x_0,z_0)$ at time $0$ of the fields $F$ and $tilde{F}$.
 \begin{lemme}
Let $(x_0,z_0)$ be an initial condition such that $x_0 \sim 0$ and $z_0 \sim 1$. For any $t$ bounded and as long as $emat^{\rho} \leq \big(\tilde{x}(t,(0,1)),\tilde{z}(t,(0,1))\big)\leq 1$ :
$$ \big(x(t,(x_0,z_0)),z(t,(x_0,z_0))\big ) \sim \big(x(t,(0,1)),z(t,(0,1))\big)$$
\end{lemme}
 \textbf{Demonstration} : We introduce the two families of vector fields defined for $\alpha >0$ by :
 
  \beq \label{eqchsignealphaplusminus}
F^{+\alpha}  \quad \left\{
\begin{array}{lcl}
\displaystyle \frac{dx}{dt}& =&1 \\[6pt]
\displaystyle \frac{dz}{dt} & =&\displaystyle (1+alpha)zf(z)
\end{array} 
\right.
\quad \quad
F^{-\alpha}  \quad \left\{
\begin{array}{lcl}
\displaystyle \frac{dx}{dt}& =&1 \\[6pt]
\displaystyle \frac{dz}{dt} & =&\displaystyle (1-\alpha)zf(z)
\end{array} 
\right.
\feq
and note $\big(x(t,(x_0,z_0),±\alpha),z(t,(x_0,z_0),±\alpha)\big)$ the solutions from initial condition $(x_0,z_0)$ at time $0$ of the vector fields $F{±\alpha}$.
On the opposite figure the segments $\arc{b,c}$, $\arc{c,d}$ and $\arc{a,e}$ which define the closed contour $\partial \Gamma_{\alpha} =\arc{a,b,c,d,e,a}$ which surrounds the domain $\Gamma_{\alpha}$ are respectively : \\
 \begin{minipage}{0.50 \textwidth} 
\bitbul
\item$\arc{b,c}$ : The portion of the trajectory of $F^{-\alpha}$ between $x = -\alpha$ and $x = 0$ that passes through the point $(0,1-\alpha)$.
\item $\arc{c,d}$ : The portion of trajectory of $F^{+\alpha}$ between $x = 0$ and $x = t$ that passes through the point $(0,1-\alpha)$.
\item $\arc{a,e}$ The portion of trajectory of $F^{-\alpha}$ comprised between $x = -\alpha $ and $x = t $ which passes through the point $(0,1+\alpha)$.
\fit

 \end{minipage}$\quad$
  \begin{minipage}{0.50 \textwidth}
 \includegraphics[width=1\textwidth]{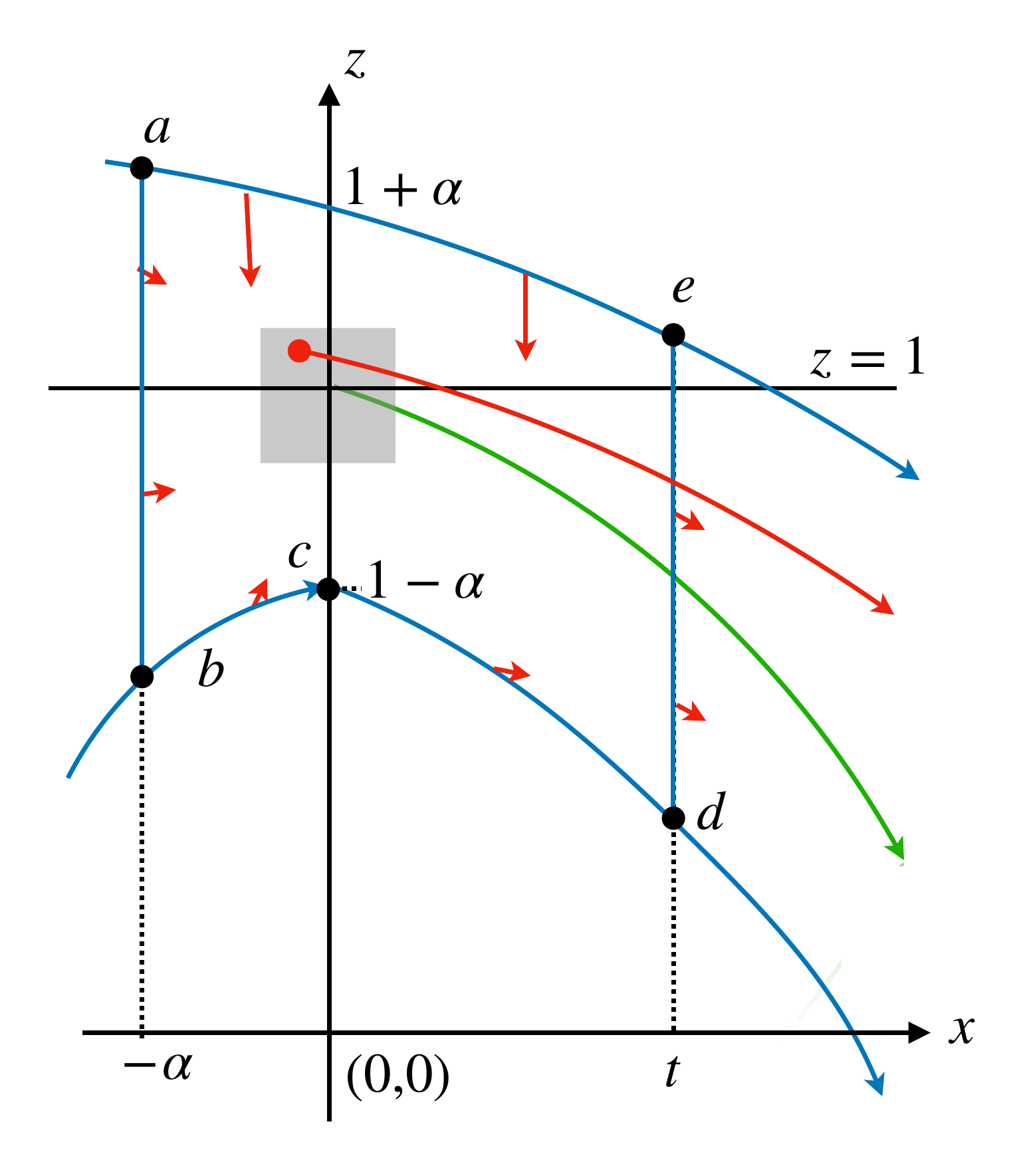}
 \end{minipage}
 The red trajectory is the solution of $F$ from an initial condition infinitely close to $(0,1)$ and the green trajectory the solution of $\tilde{F}$ from $(0,1)$.
We can easily verify that if $\alpha \gnsim 0$ the field $F$ points strictly towards the interior of $\Gamma_{\alpha}$ on all the segments of the contour $\partial \Gamma_{\alpha}$ except on $\arc{d,e}$ ; on the other hand, always if $\alpha \gnsim 0$ the point $(x_0,z_0) \sim (0,1)$ belongs to $\Gamma_{\alpha}$ and therefore (c. f. the theorem \ref{piege}) the trajectory of $F$ which comes from it 
exits $\Gamma_{\alpha}$ at a point on the segment $\arc{d,e}$. As the fields $F^{-\alpha}$ and $F^{+\alpha}$ tend to $F$ and, given the 
assumptions of regularity on $f$, when $\alpha$ tends to $0$ the length of the segment $\arc{d,e}$ tends to $0$ and consequently $ |z(t,(x_0,z_0))-\tilde{z}(t,(0,1))|\leq \delta$ for all $$\delta \gnsim 0$$ which completes the demonstration.\\
$\Box$


\fin

